\documentclass[12pt]{article}
\usepackage[round,sort,authoryear]{natbib}
\usepackage{cancel}
\usepackage{amsmath,amssymb,amsthm,amsfonts, fullpage, setspace}
\usepackage{eucal} 
\usepackage{subeqnarray,bm}
\usepackage[utf8]{inputenc}
\usepackage{graphicx}[final]
\usepackage{psfrag} 
\usepackage{color,float}
\usepackage{pdfpages}
\usepackage{epstopdf}
\usepackage[colorlinks]{hyperref}
\usepackage{hyperref}
\hypersetup{colorlinks=red,citecolor=blue,pdfstartview=FitH, pdfpagemode=None}
\usepackage{multirow}
\usepackage[bottom]{footmisc}
\usepackage{ifthen} 
\usepackage{sectsty}
\usepackage{subfig}
\usepackage{float}
\usepackage[english]{babel}
\usepackage[nottoc]{tocbibind}
\usepackage{calrsfs}
\usepackage{soul}
 \DeclareGraphicsExtensions{.eps,.png,.pdf}
\usepackage{mwe}
\makeatletter
\def\maxwidth{ %
  \ifdim\Gin@nat@width>\linewidth
    \linewidth
  \else
    \Gin@nat@width
  \fi
}
\makeatother

\definecolor{fgcolor}{rgb}{0.345, 0.345, 0.345}

\allsectionsfont{\normalsize\bfseries}

\usepackage{framed}
\makeatletter
 {\par\unskip\endMakeFramed%
 \at@end@of@kframe}
\makeatother

\definecolor{shadecolor}{rgb}{.97, .97, .97}
\definecolor{messagecolor}{rgb}{0, 0, 0}
\definecolor{warningcolor}{rgb}{1, 0, 1}
\definecolor{errorcolor}{rgb}{1, 0, 0}

\newtheorem{thm}{Theorem}[section]

\newtheorem{defn}[thm]{Definition}

\newtheorem{remark}[thm]{Remark}

\newcommand{\ds}{\displaystyle}

\newcommand{\norm}[1]{\left\Vert#1\right\Vert}
\newcommand{\abs}[1]{\left\vert#1\right\vert}
\newcommand{\set}[1]{\left\{#1\right\}}

\newcommand{\R}{\mathbb{R}}

\newcommand{\Z}{\mathbb{Z}}
\newcommand{\p}{\partial}
\newcommand{\Pp}{\mathbb{P}}

\numberwithin{equation}{section}

\usepackage{calrsfs}

\allsectionsfont{\normalsize\bfseries}

\numberwithin{equation}{section}
\allowdisplaybreaks

\usepackage[table, dvipsnames]{xcolor}
\DeclareMathAlphabet{\pazocal}{OMS}{zplm}{m}{n}
    \newcounter{example}[section]

\usepackage{authblk}
\title{Wave Propagation  Dynamics via Lattice Difference Equations }
\author{Eddy Kwessi\footnote{ ekwessi@trinity.edu, Trinity University, Department of Mathematics, San Antonio, Texas 78212}}
\date{}	
\begin{document}
\maketitle
\begin{abstract}
We develop and analyze a lattice difference equation (LDE) framework to model the spatial dynamics of invasion in populations. This framework extends beyond classical integro-difference and reaction-diffusion models by incorporating spatial discreteness and habitat fragmentation more faithfully, making it well-suited for urban and patchy landscapes.  We characterize the local stability of equilibria, and demonstrate the existence of traveling wave solutions. A key focus is on how dispersal kernels--ranging from Gaussian to Cauchy--interact with the Allee effect to influence wave formation, propagation speed, and invasion success.  Our numerical simulations reveal that long-tailed kernels can overcome the Allee threshold through seeding effects, significantly accelerating wave fronts. These findings have direct implications for vector control strategies, informing optimal release thresholds and spatial targeting in heterogeneous environments. Furthermore, we derive a stochastic characterization of outbreak size as a Poisson-binomial distribution, offering probabilistic insight into local infection burden and paving the way to a novel minimization criteria of release strategies based on bi-modality.  Our results provide a mathematically grounded basis for predicting spatial spread and optimizing release strategies for {\it Wolbachia}-based vector control programs. 
This work bridges theory and application, providing both analytical insights and computational tools for understanding spatial epidemiology in discrete habitats
\end{abstract}

{\bf Key words}: Lattice Difference Equations, {\it Wolbachia}, Allee Effect, Traveling Waves, Spatial Epidemics, Outbreak Distribution.

\section{Introduction}\label{sec1}

{\it Wolbachia} is an endosymbiotic bacterium that infects a wide range of arthropods and has gained prominence as a biological control agent for vector-borne diseases such as dengue, Zika, chikungunya, and malaria \citep{hoffmann2011successful, werren2008wolbachia}. By manipulating host reproductive mechanisms—most notably through cytoplasmic incompatibility (CI)—and by reducing the competence of mosquitoes to transmit pathogens, {\it Wolbachia} has emerged as a cost-effective and sustainable intervention strategy \citep{ferguson2015modeling, turelli2010cytoplasmic}. Recent field trials have demonstrated its potential for suppressing disease transmission in urban and semi-urban settings \citep{oneill2018scaled}.

Understanding the spatial dynamics of {\it Wolbachia} invasion is critical for optimizing its deployment and predicting long-term stability. Classical ordinary differential equation (ODE) models have offered foundational insights into {\it Wolbachia} infection dynamics based on maternal transmission and fitness trade-offs \citep{dobson2002wolbachia, hancock2011modeling}. However, these models lack spatial resolution and cannot capture heterogeneities in landscape structure or dispersal patterns, which play a crucial role in the success or failure of biological invasions.

Spatial models--particularly reaction-diffusion equations and integro-difference equations (IDEs)--have addressed this gap by incorporating mosquito dispersal and environmental variability. IDEs are especially suited to model mosquito populations due to their discrete generational structure and ability to simulate nonlocal interactions through dispersal kernels \citep{neubert2000demography}.

However, IDEs may still oversimplify population structure in fragmented or patchy environments, where mosquito habitats are inherently discrete. To address this limitation, we introduce a lattice difference equation (LDE) framework for modeling the spatial spread of {\it Wolbachia}. 

LDEs represent spatial domains as a grid or lattice, enabling the modeling of mosquito movement between discrete locations and allowing for a more faithful representation of local heterogeneity, habitat fragmentation, and non-uniform dispersal \citep{weinberger1982long, li2020lattice}.

This discrete spatial structure is particularly relevant for modeling invasion fronts and understanding thresholds in bi-stable systems driven by Allee effects \citep{skellam1951random}. We also incorporate a biologically motivated {\it Wolbachia} growth function derived from CI intensity and relative fitness costs, which supports a bi-stable dynamic characterized by extinction, threshold, and persistence equilibria.

The incorporation of diverse dispersal kernels (e.g., Gaussian, Cauchy, power-law) allows us to evaluate how local versus long-range dispersal affects wave propagation, and invasion success \citep{fife1979mathematical}.

This paper aims to rigorously study the dynamics of spread in lattice-based mosquito populations by:
\begin{enumerate}
    \item Analyzing local stability near biologically meaningful fixed points;
    \item Proving the existence of traveling waves and wavefront propagation;
    \item Investigating how dispersal mechanisms interact with Allee effects to accelerate or inhibit invasion;
    \item Exploring the distribution of outbreak size and using its modality shift as a novel criterion for release cost minimization.
\end{enumerate}

By developing this lattice-based framework, we offer new insight into the spatial epidemiology of {\it Wolbachia} and provide tools for designing efficient release strategies in heterogeneous environments.

%%
%%%%%%%%%%%%%%%%%%%%%%%%%%%%%%%%%%%%%%%%%%%%%%%%%%%%%%%%%%%%%%%
\section{Preliminaries}\label{sec:prelim}
%%%%%%%%%%%%%%%%%%%%%%%%%%%%%%%%%%%%%%%%%%%%%%%%%%%%%%%%%%%%%%%

In this paper, we will be chiefly concerned with  the {\it Wolbachia}  growth function for (see \cite{Yu2019})  is given  as 
\begin{equation}\label{eqn0}
f(v(t))=\frac{(1-s_f)v(t)}{s_hv^2(t)-(s_h+s_f)v(t)+1}\;,
\end{equation}
where $s_h$ represents the cytoplasmic incompatibility (CI) intensity,  $s_f$ represents the relative fitness cost of infected females/males, and $v(t)$ represents the  infection frequency at time $t$. 
%However, for comparison purposes, we will  sometime discuss the above growth function conjointly with  the Allee growth function given as 
%\begin{equation}\label{eqn00}
%f(v(t))=rv((t)\left(1-\frac{v(t)}{K}\right)\left(\frac{v(t)}{A}-1\right)\;,
%\end{equation}
%where  $r$ represents the intrinsic growth rate, $K$ represents the carrying capacity, and   $A$ represents the Allee threshold.
Integro-difference equations (IDEs) provide a powerful framework for modeling the spatial spread of {\it Wolbachia} by incorporating both local reproductive dynamics and long-range dispersal. These models account for the fact that mosquito populations undergo generational turnover with discrete time steps $t$, making them well-suited for studying wavefront propagation in heterogeneous environments. Specifically, the {\it Wolbachia} spread model can be extended to a spatially explicit setting using an integro-difference equation of the form (see  \cite{kot1996dispersal}):
    \begin{equation}\label{eqn1}
    v(t+1,x) = \int_{-\infty}^{\infty} K(x, y) f(y,v(t,y)) dy,
    \end{equation}
where $v(t,y)$  represents the infection frequency at location  $y$ and time or generation $t$, $K(x,y)=K(x-y)$ is the dispersal kernel capturing mosquitoes movements (relative distance from $x$ to $y$). More specifically, it is  the probability of moving from location $y$ to $x$. The function $f(x,v(t,x))$ represents the local growth function determined by cytoplasmic incompatibility and fitness costs. 
If only a fraction $\delta$ of  the population disperses from location $y$ to location $x$, the above could be extended to  obtain the integro-difference equation:
  \begin{equation}\label{eqn2}
    v(t+1,x) = (1-\delta(x))f(x,v(t,x))+\int_{-\infty}^{\infty} \delta(y) K(x, y) f(y,v(t,y)) dy\;.
    \end{equation}
Despite the usefulness of IDEs, they may not fully capture population dynamics in environments where mosquito habitats are highly fragmented, leading to naturally discrete population distributions. An alternative approach using lattice difference equations (LDE) allows for a more granular, spatially structured model. The lattice difference equation framework considers a discrete set of spatial locations  and models the spatial spread. Now, let $d$ be a positive. We define integer the $d$-dimensional integer lattice as  
 \[
\Z^d=\set{(k_1,k_2,\cdots, k_d): k_\ell \in \Z, \ell \in \set{1, 2, 3,\cdots, d}}\;.\]
\begin{figure}[H]
 \resizebox{1\textwidth}{!}{\begin{minipage}{1\textwidth}
\centering \begin{tabular}{cc}
\includegraphics[scale=0.35]{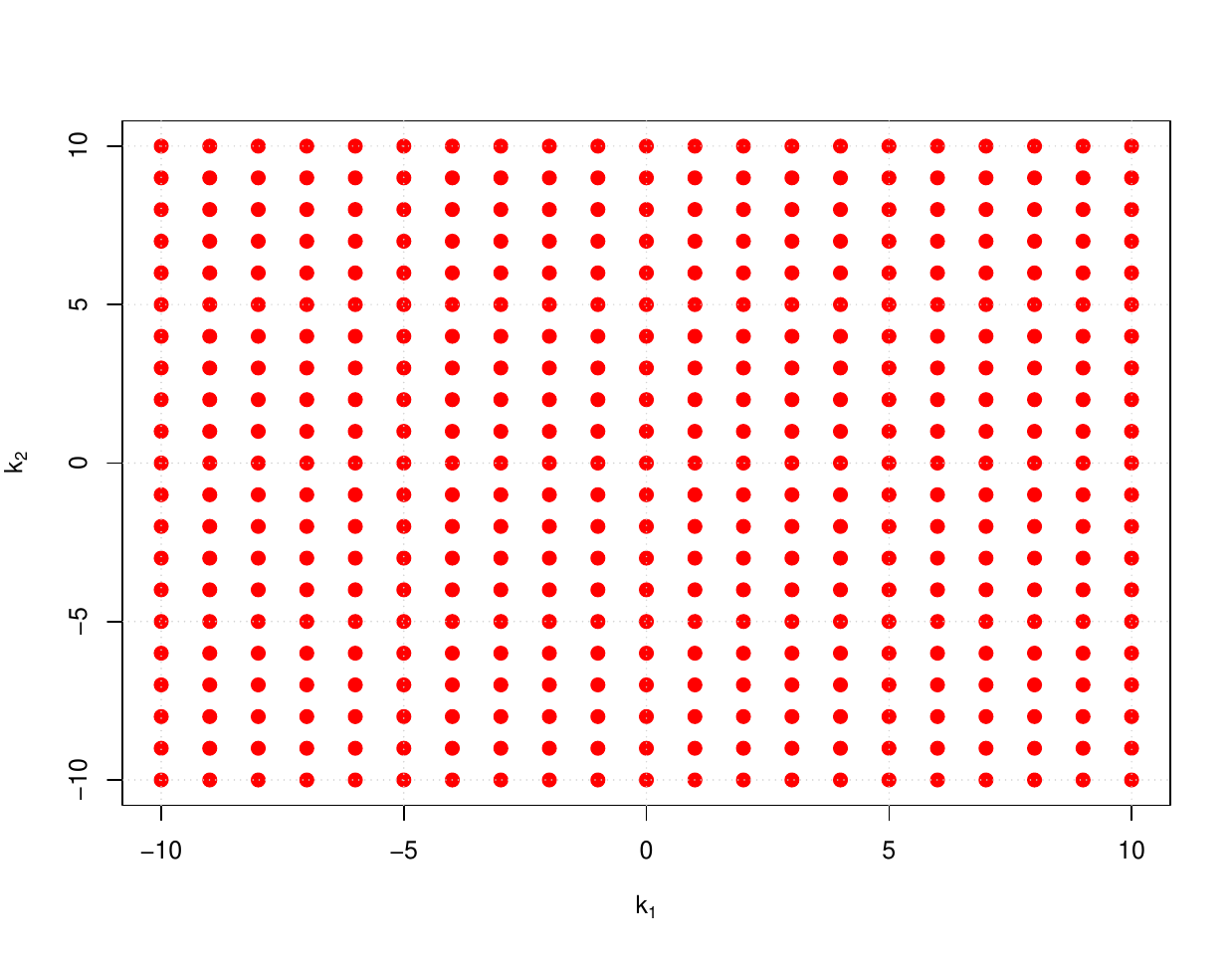} & \includegraphics[scale=0.25]{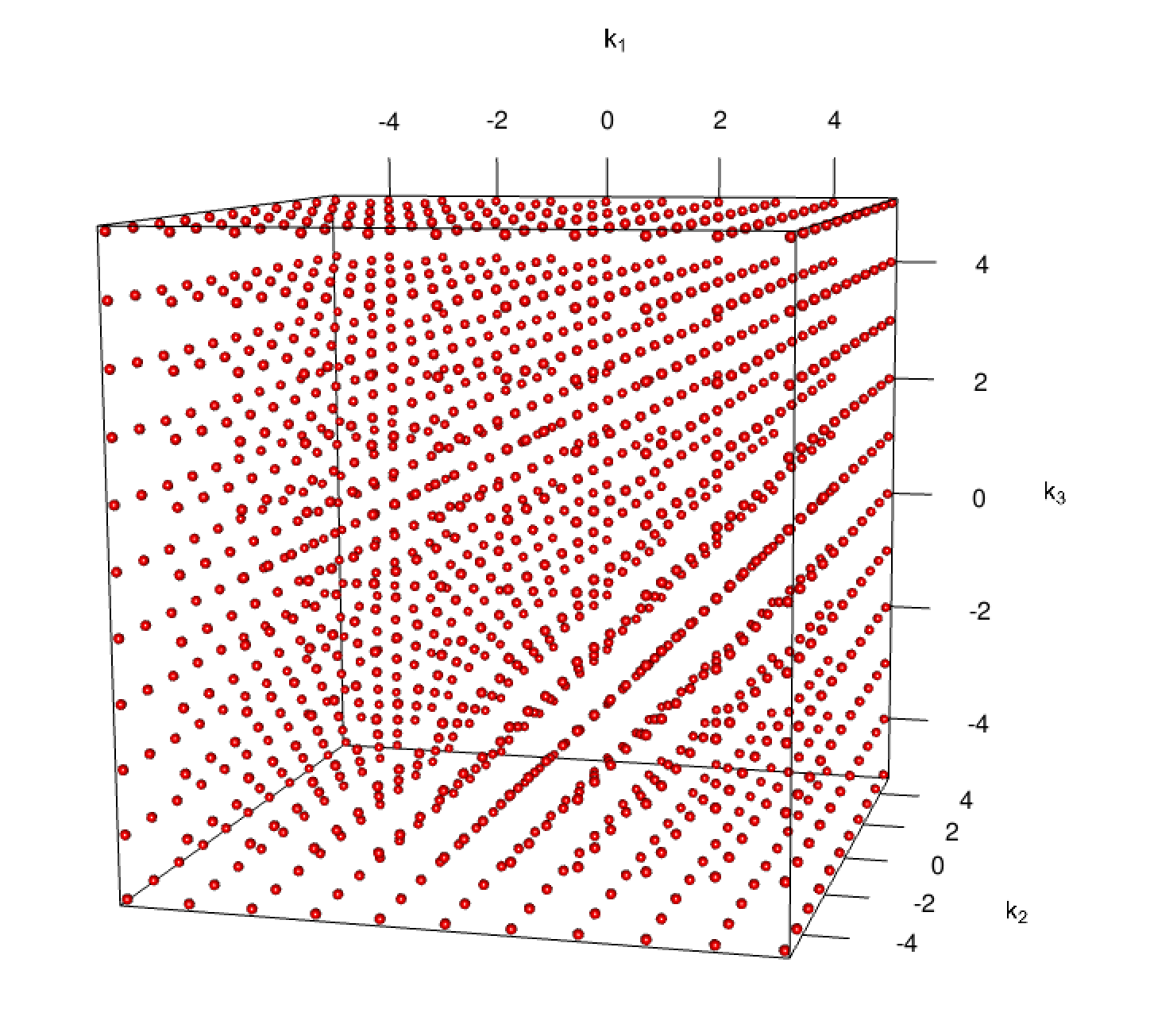}
\end{tabular}
\end{minipage}}
\caption{Representations of lattices for $d=2$ (left) and $d=3$ (right).}
\label{fig:1}
\end{figure}

 \noindent For $i,j\in \Z^d$, put
 \begin{eqnarray*}
 \delta_i&=&\delta(x_i)\;,\\
  K_{ij}&=& K(x_i,y_j)\;,\\
 v_i(t)&=&v(t,x_i)\;,\\
 f_i(v_i(t))&=&f(x_i,v(t,x_i))\;.
 \end{eqnarray*}
We will consider the autonomous Wolbachia lattice difference equation (WLDE), obtained from the spatial discretization of \eqref{eqn2} above, given as 
 \begin{equation}\label{eqn3}
 v_i(t+1)=(1-\delta_i)f_i(v_i(t))+\sum_{j\in \Z^d}\delta_jK_{ij}f_j(v_j(t)),\quad \mbox{for $i\in \Z^d$}\;,
 \end{equation}
where $\delta_i$ controls the proportion of the population that disperses at location $x_i$, and $K_{ij}$ defines the movement probabilities between sites $x_i$ and $x_j$. We observe that $\Z^d$ defines a spatial grid on which   each node represents a discrete location where mosquito populations reside. This formulation allows for a more accurate representation of invasion dynamics in spatially structured mosquito populations.

%%%%%%%%%%%%%%%%%%%%%%%%%%%%%%%%%%%%%%%%%%%%%%%%%%%%%%%%%%%%%%
%%%%%%%%%%%%%%%%%%%%%%%%%%%%%%%%%%%%%%%%%%%%%%%%%%%%%%%%%%%%%%
\section{Stability Analysis}\label{sec:stability analysis}
%%%%%%%%%%%%%%%%%%%%%%%%%%%%%%%%%%%%%%%%%%%%%%%%%%%%%%%%%%%%%%%
%%%%%%%%%%%%%%%%%%%%%%%%%%%%%%%%%%%%%%%%%%%%%%%%%%%%%%%%%%%%%%
Before we discuss the stability of the above LDE, let us discuss its fixed points. The fixed points of the above LDE satisfy $v_i(t+1)=v_i(t)$ for all $i\in \Z^d$. This amounts to 
\[
v_i(t)=(1-\delta_i)f_i(v_i(t))+\sum_{j\in \Z^d}\delta_jK_{ij}f_j(v_j(t))\;.
\]
For $\delta_i=\delta=$constant, we have
\[
f(v_i(t))-v_i(t)=\delta\left(f(v_i(t))-\sum_{j\in \Z^d}K_{ij}f_j(v_j(t))\right)\;.
\]
We see that if $v_i(t)$ is a fixed point of $f$, this further reduces to 
\[
v_i(t)=\sum_{j\in \Z^d}K_{ij}v_j(t)\;.
\]
In vector form, this amounts to $Kv=v$, where $K=(K_{ij})_{i,j\in \Z^d}$.
The latter equation means that a fixed point $v(t)=(v_i(t))_{i\in \Z^d}$ for the LDE is a fixed point of the Wolbachia growth function  and  the  eigenvector of the stochastic matrix $K$ associated the eigenvalue 1. We know by the Perron-Frobrenius Theorem that such a vector exists. Examining the growth function $f$, its fixed points are respectively $v^*=0, v^*=\frac{s_f}{s_h}$ and $v^*=1$. 
%The  fixed point $v^*=\frac{s_f}{s_h}$ is referred to as the Allee threshold. 
Now we state the  result below about the stability of the LDE around a fixed point $v^*$.
\begin{thm}\label{thm:persist}
Let \( \widehat{K}\) be the  Fourier transform  of the dispersal kernel \( K \). Then we have the following:
\begin{enumerate}
\item If \(\ds \sup_{k\in (\Z^+)^d} |f'(v^*)| |(1 - \delta) + \delta  \widehat{K}(k) | < 1\), then the fixed point $v=(v^*_k)$ is locally asymptotically stable (LAS).
\item If \(\ds \sup_{k\in (\Z^+)^d} |f'(v^*)| |(1 - \delta) + \delta \widehat{K}(k) | >1\), then fixed point $v=(v^*_k)$ is  unstable (UNS).
\end{enumerate}
\end{thm}
The proof of this theorem can be found in Appendix A.

\subsection{Dispersal Kernel Choices}
\begin{figure}[H] %  figure placement: here, top, bottom, or page
   \centering
   \includegraphics[scale=0.5]{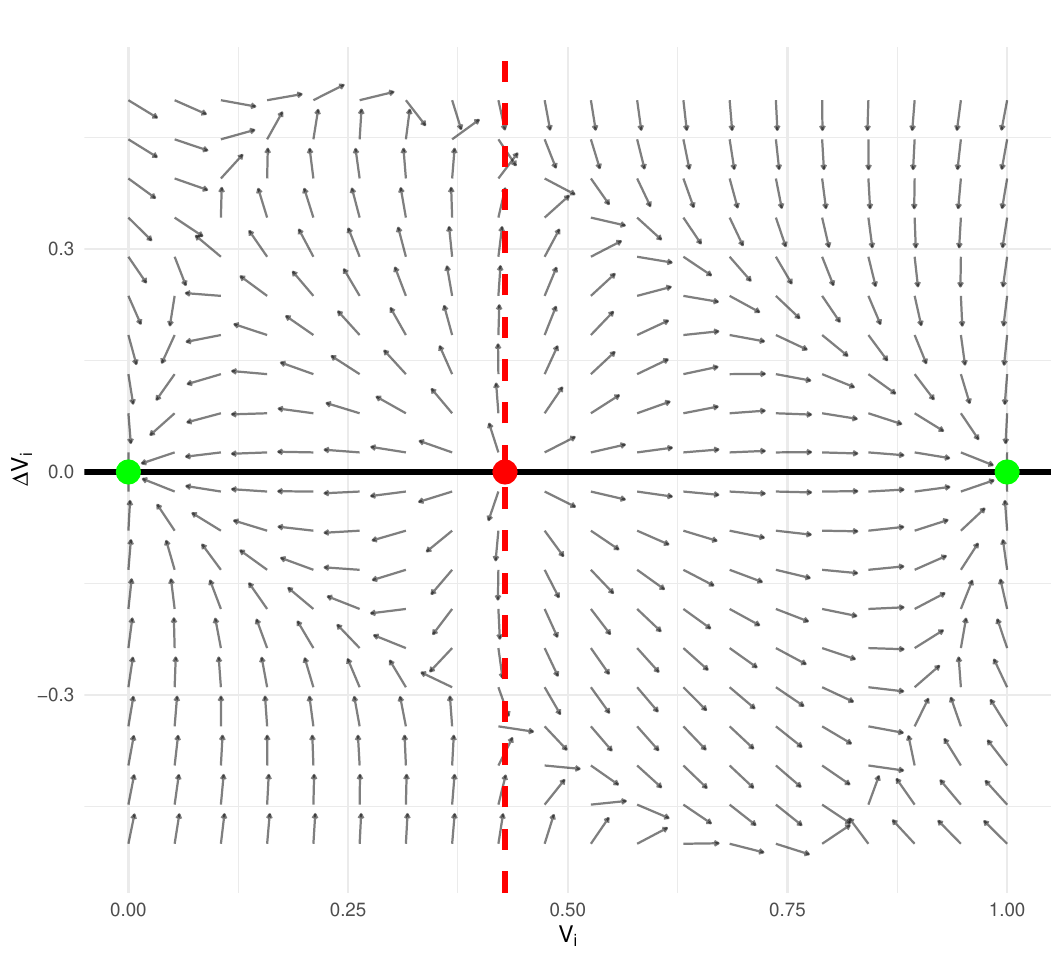} 
   \caption{Phase space diagram in the space $(V_i,\Delta V_i)$, obtained for $s_f=0.3, s_h=0.7$, and $\delta=0.6$. The green dots represents respectively $v6*=0$ and $v^*=1$, while the red dot represent the Allee Threshold $v^*=\frac{s_f}{s_h}$. The arrows represents the vector field showing how the infection progresses in the space.}
   \label{fig:PSD}
\end{figure}

In practice, one may  use one of the  kernels in Table \ref{tab1} below, with their corresponding Fourier Transform for $d=2$.
\begin{table}[h]
\resizebox{1.15\textwidth}{!}{\begin{minipage}{1\textwidth}
\begin{tabular}{|l|l|l|}
\hline
{\bf Name} &Kernel & DFT\\ \hline
{\bf Cauchy} &$\ds C(m, n) = \frac{\gamma}{\pi (\gamma^2 + m^2 + n^2)}\;, \mbox{for some $\gamma >0$}$& $\ds \widehat{C}(k, l) = \frac{\exp \left( -2\pi \gamma \sqrt{k^2 + l^2} \right)}{1 + e^{-2\pi \gamma}}$\\ \hline
{\bf Power Law}&$\ds P(m, n) = \frac{C}{(1 + m^2 + n^2)^{\gamma/2}}$& $\ds \widehat{P}(k, l)  \approx   \frac{C}{(1 + \frac{k^2}{M^2} + \frac{l^2}{N^2})^{\gamma/2}} $\\ \hline
{\bf Gaussian} &$\ds G(m, n) = \frac{1}{2\pi \sigma^2} \exp\left(-\frac{m^2 + n^2}{2\sigma^2} \right)$&$\ds  \widehat{G}(k, l) = \exp \left( -2\pi^2 \sigma^2 \left( \frac{k^2}{M^2} + \frac{l^2}{N^2} \right) \right)$\\ \hline

{\bf Uniform} &$\ds U(m, n) = \frac{1}{MN}, \quad \text{for } 0 \leq m < M, \, 0 \leq n < N$& $\ds \widehat{U}(k, l) =\begin{cases}
\ds  \frac{1}{MN} \left( e^{-\pi \hat{j} \frac{k}{M}}e^{\pi \hat{j} \frac{l}{N}}\frac{\sin{\pi k}}{\sin{\frac{\pi k}{M}}} \cdot 
\frac{\sin{\pi l}}{\sin{\frac{\pi l}{N}}} \right) & \textrm{if $k\neq 0\neq l$}\\
1 & \textrm{if $k=0=l$}
\end{cases}$\\ \hline
\end{tabular}
\end{minipage}}
\caption{Commonly used kernels and their respective Fourier Transform, where $\hat{j}^2=-1$.}
\label{tab1}
\end{table}
The Allee effect and dispersal mechanisms jointly influence the speed and behavior of  waves in Wolbachia invasion models. Below, we examine how these two factors interact to either accelerate or hinder wave propagation.  We use Theorem \ref{thm:persist} above to summarize  the dynamics of the  WLDE, for kernels such that $\ds \sup_{k\in (Z^+)^d}~\widehat{K}(k)=1$.\\
At $v^*=0$ we have that $|f'(0)|=|1-s_f|<1$ which implies that local asymptotic stability. At $v^*=1$, we have that $|f'(1)|=\left|\frac{1-s_f}{1-s_h}\right|<1$ which also implies that local asymptotic stability.  At $v^*=\frac{s_f}{s_h}$, we have  $\ds f'\left(\frac{s_f}{s_h}\right)=\frac{s_h-s_f^2}{s_h-s_hs_f}>1$ which implies that the point in unstable. In fact, we see  in Figure \ref{fig:PSD} below representing the  phase space diagram  that this point is in fact a repeller.
This phase space diagram clearly illustrates the dynamics of the WLDE around the fixed points. Noteworthy is the fact that for $\Delta V_i$ large enough, it is possible to connect the two stable fixed points $v^*=0$ and $v^*=1$, jumping the unstable fixed point $v^*=\frac{s_f}{s_h}$. Orbits connecting these two fixed points are called traveling waves and their existence will be discussed in the section below.

%%%%%%%%%%%%%%%%%%%%%%%%%%%%%%%%%%%%%%%%%%%%%%%%%%%%%
\section{Traveling waves} \label{sec:travelling waves} 
%%%%%%%%%%%%%%%%%%%%%%%%%%%%%%%%%%%%%%%%%%%%%%%%%%%%%

Traveling waves represent a spatially structured invasion where the {\it Wolbachia}-infected population spreads steadily across a habitat. These waves determine how fast and under what conditions the infection propagates, which has biological and epidemiological significance in disease control. We will make the following assumption of the kernel
\[A_1: \mbox{For all $i\in \Z^d$, we have  $\ds \sum_{i\in \Z^d}K_{ij}=1$}\;.\]

%%%%%%%%%%%%%%%%%%%%%%%%%%%%%%%%%%%%%%%%%%%%%%%%%%%%%%%%%%%%%%
\subsection{Existence of Traveling Waves}
%%%%%%%%%%%%%%%%%%%%%%%%%%%%%%%%%%%%%%%%%%%%%%%%%%%%%%%%%%%%%%%
In this section, we will discuss the existence of traveling waves for the WLDE given in equation \eqref{eqn3} with constant dispersal probability, that is, $\delta_i=\delta=\mbox{constant}$. We first observe that since $0<\delta\leq 1$ and $0\leq f(v_i(t))\leq 1$ for all $i\in \Z^d$, then equation \eqref{eqn3} is well-defined for dispersal kernels $K$ satisfying $A_1$.
We begin with some important definitions:
\begin{defn} A traveling wave for {\it Wolbachia}-LDE is a solution written using the ansatz
\[v_i(t)=U(\xi),\quad \xi=i-ct\;,\]
where 
\begin{itemize}
\item $U(\xi)$ is called the wave profile,
\item $\xi=i-ct$ is the traveling wave coordinate,
\item $c$ is the wave speed. 
\end{itemize}
\end{defn}
\noindent Substituting $v_i(t)$ into the equation, we obtain the equation of the traveling wave equation 
\begin{equation}\label{eqn8}
U(\xi+c)=(1-\delta)f(U(\xi))+\delta\sum_{j\in Z^d} K_{ij}f(U(\xi+i-j))\;.
\end{equation}
This equation expresses how the wave propagates forward in the spatial domain. In biological applications, traveling waves often represent the invasion of {\it Wolbachia}-infected mosquitoes into an uninfected population. Thus, we impose biologically meaningful boundary conditions to ensure that the {\it Wolbachia} strain invades and persists.\\
%\noindent {\bf Assumptions}\\
\begin{itemize}
\item[$A_2$:] $\ds \lim_{\xi\to \infty}U(\xi)=0$. This means that the leading edge of the wave $\xi\to \infty$ corresponds to the uninfected mosquito population.
\item[$A_3$:] $\ds \lim_{\xi\to -\infty}U(\xi)=1$. This means that the trailing  edge of the wave $\xi\to -\infty$  corresponds to a fully established {\it Wolbachia}-infected population.
\end{itemize}
We can now state the result of on the existence of traveling waves. 
\begin{thm} \label{thm:1} If Assumptions $A_1-A_3$ are met, then the {\it Wolbachia}-LDE in \eqref{eqn3} possesses traveling waves.
\end{thm}
\begin{figure}[H] %  figure placement: here, top, bottom, or page
   \centering
   \includegraphics[scale=0.7]{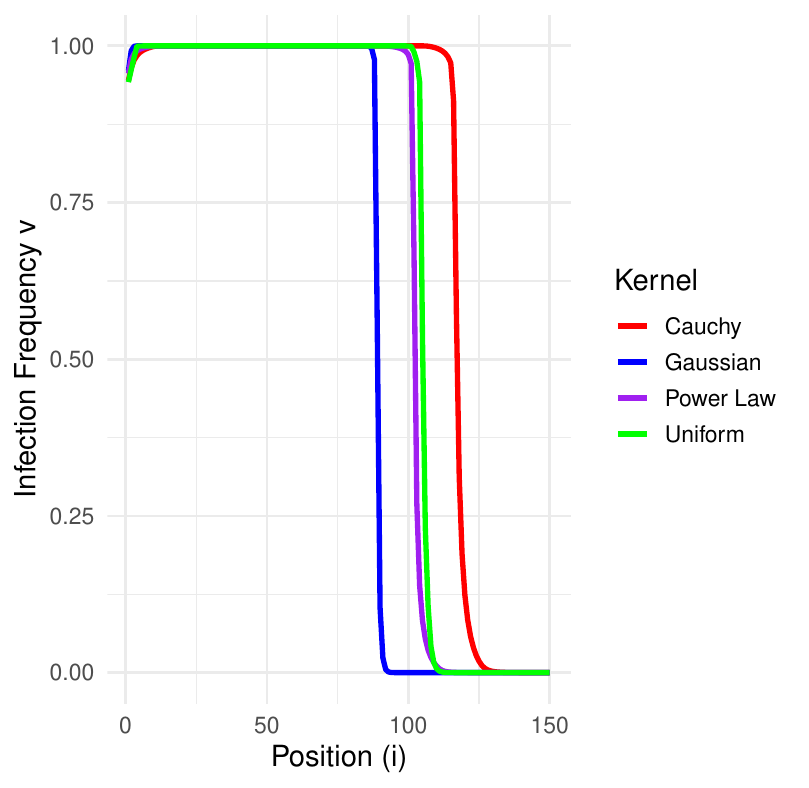} 
   \caption{This figure was obtained on 1 dimensional lattice with width 150 and for $sh=0.8, A=0.4,\delta=0.1, $ and over a discrete times $t=0,1,\cdots, 200$, evaluated at time $t=150$. It shows traveling waves for different kernel types, connecting the two stable fixed points $v^*=1$ and $v^*=0$.  Although the shape of the traveling waves is the same for all four kernels,  the Cauchy kernel seems to be the one with the farthest reach. We also observe that tuning these parameters may produce different outcomes. }
   \label{fig:travelWave}
\end{figure}
\begin{remark} 1. It is important to note that LDE under consideration here has a stable equilibrium at $v=0$. Because the growth function $f(v)$ possesses an  Allee effect, this implies that small perturbations around $v=0$ decay, so the population/infection does not invade if it starts below the Allee threshold $v^*$. Therefore, waves only form and propagate if the initial condition exceeds the threshold somewhere. The wave speed is not determined by the behavior at 
$v=0$ because disturbances near 0 shrink rather than grow.  This means that the wave is ``pushed" from behind, not pulled by growth at the front. Consequently, the technique of linear determinacy often used to determine the minimal wave speed does not apply here. Hence the minimal traveling wave speed may  be determined using methods such as the phase plane analysis, numerical computation or constructive wave ansatz.\\
2. In a pure difference equation or local map $V_{t+1}=f(v_t)$, the flow is one-dimensional and monotonic, that is, trajectories cannot jump over an unstable fixed point. The Allee threshold $A$ blocks spread from infinitesimal introduction. However in spatial system like the Wolbachia. If  the invasion  exceeds $A$ over a sufficiently wide region, the wave can self-propagate, creating a heteroclinic connection in function space. 
\end{remark}
\subsection{Wave speed}
An important aspect in LDE is the wave speed $c(t)=\frac{dx}{dt}$. The wave speed  is critically important because it quantifies the rate at which the Wolbachia infection spreads through the mosquito population across space. To fix ideas, we recall that if $v=v(t,x)$, then \[\frac{dv}{dt}=\frac{\partial v}{\partial x}\frac{dx}{dt}+\frac{\partial v}{\partial t}\frac{dt}{dt}=\frac{\partial v}{\partial x}\frac{dx}{dt}+\frac{\partial v}{\partial t}.\] 
We have \[\frac{dv}{dt}=0\implies c(t)=\frac{dx}{dt}=-\frac{\frac{\partial v}{\partial t}}{\frac{\partial v}{\partial x}}\;.\]
In the context of difference equations, we can make the following approximations: 
\[\frac{\p v}{\p t}\approx v(x,t+1)-v(x,t) \]
And 
\[\frac{\p v}{\p x}\approx v(x+1,t)-v(x,t) \;.\]
It follows that if $v(x+1,t)-v(x,t) \neq 0$, we define  the wave speed at time $t$ as: 
\[ c(t):=-\frac{ v(x,t+1)-v(t,x) }{v(x+1,t)-v(x,t)}\;,\]
and the asymptotic wave speed as: 
\[c^*=\lim_{t\to \infty}c(t)\;.\]
We observe that $v(x, t+1)=v(x,t)$ if $v(x,t)=v^*=v(x^*(t),t)$, a fixed point of the LDE in \eqref{eqn2}. In this case, we clearly have $c(t)=0$. 
Given a kernel, the key parameters of the LDE are the Allee threshold and the dispersion rate $\delta$. Let us illustrate this using the Wolbachia growth function. Of interest are: the effect of the local dispersion parameter $\delta$, the effect of the Allee Threshold, and the effect of the initial wave speed on the asymptotic wave speed. 
\subsubsection{Effect of the local dispersion parameter}
\begin{figure}[H] %  figure placement: here, top, bottom, or page
   \centering
   \includegraphics[scale=0.8]{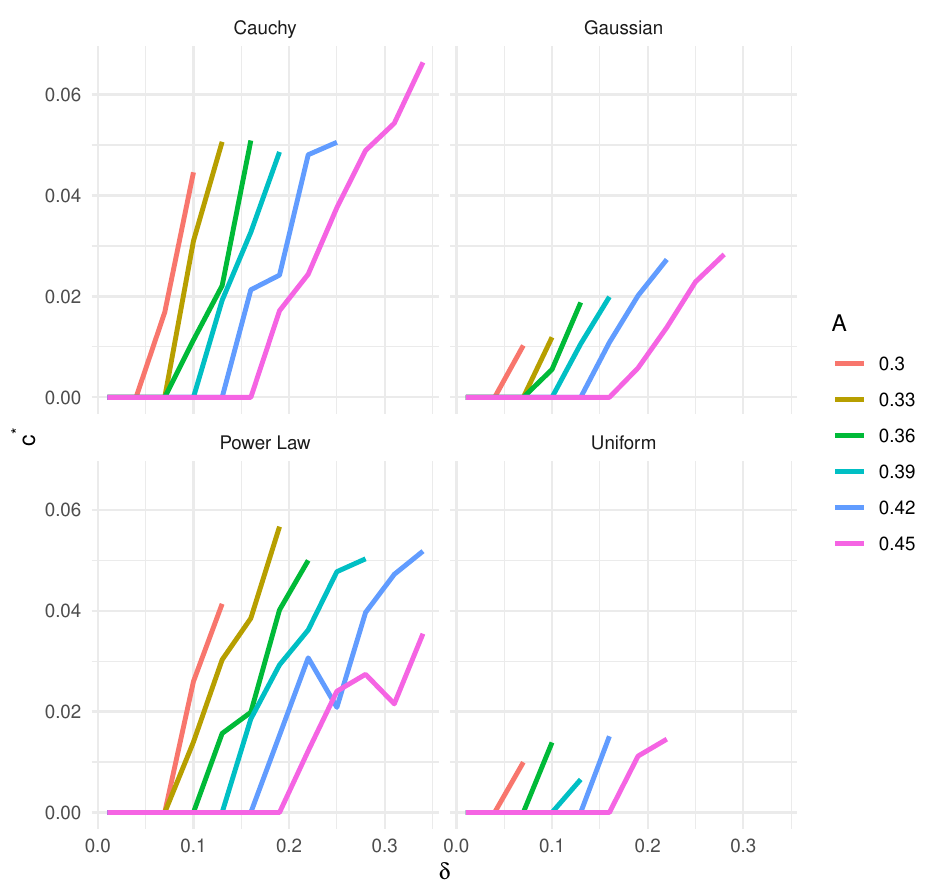} 
   \caption{This plot shows how  the asymptotic wave speed is affected for selected values  of the Allee threshold $A$ and values of $\delta \in [0.01,0.3]$. We observe that as the local dispersion $\delta$ increases, so does the asymptotic speed. Moreover, the Cauchy dispersal kernel  tends to produce higher speed than the other three kernels. }
   \label{fig:AsymptoticSpeedVsDelta}
\end{figure}

For fixed $A=0.4$ and $s_h=0.8$, we plot $c^*$ as function of $\delta$, see Figure  \ref{fig:AsymptoticSpeedVsDelta}. 

\subsubsection{Effect of the Allee threshold}
\begin{figure}[H] %  figure placement: here, top, bottom, or page
   \centering
   \includegraphics[scale=0.8]{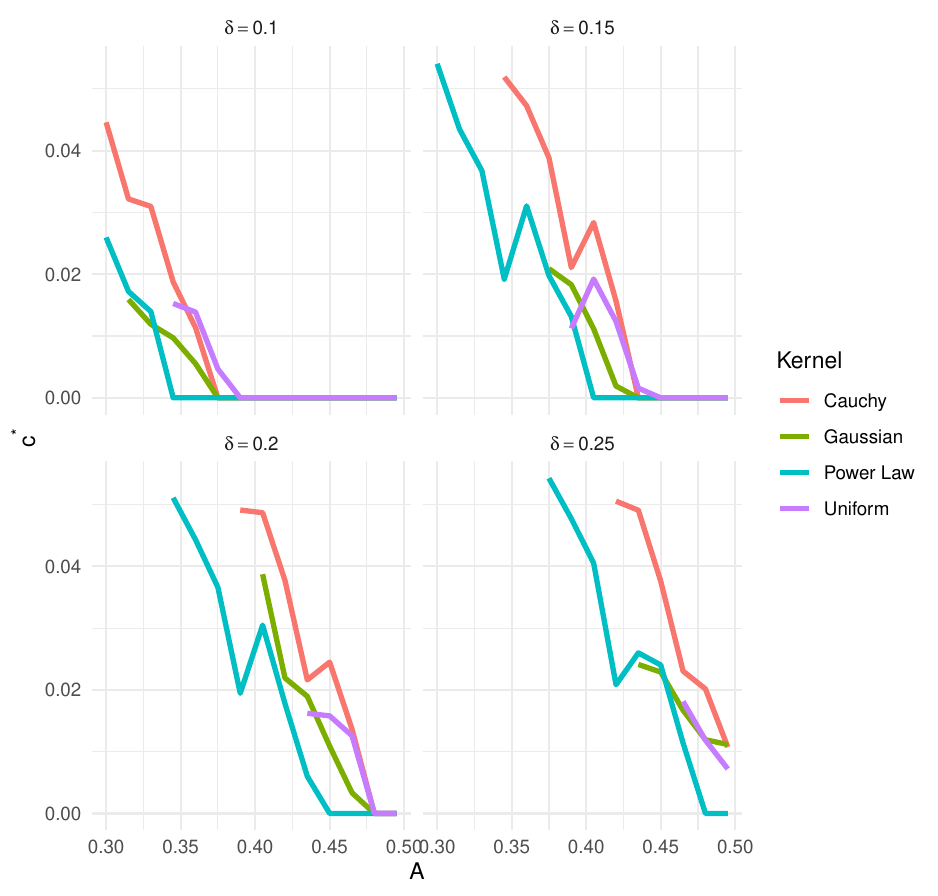} 
   \caption{These figures show the effect of the Allee threshold for different values of the local dispersion on the asymptotic speed. Large values of $A$ ultimately lead the waves to disappear.  The Cauchy kernel, as a long-range dispersal kernel is superior to all others because its asymptotic related speeds are much higher and can handle much higher Allee effect thresholds.}
   \label{fig:AsymptoticSpeedVsA}
\end{figure}

For fixed $\delta=0.3$, we plot the asymptotic waves speed $c^*$ as a function of $A$, see Figure \ref{fig:AsymptoticSpeedVsA} . 
\subsubsection{Effect of the initial release speed}
\begin{figure}[H] %  figure placement: here, top, bottom, or page
   \centering
   \includegraphics[scale=0.6]{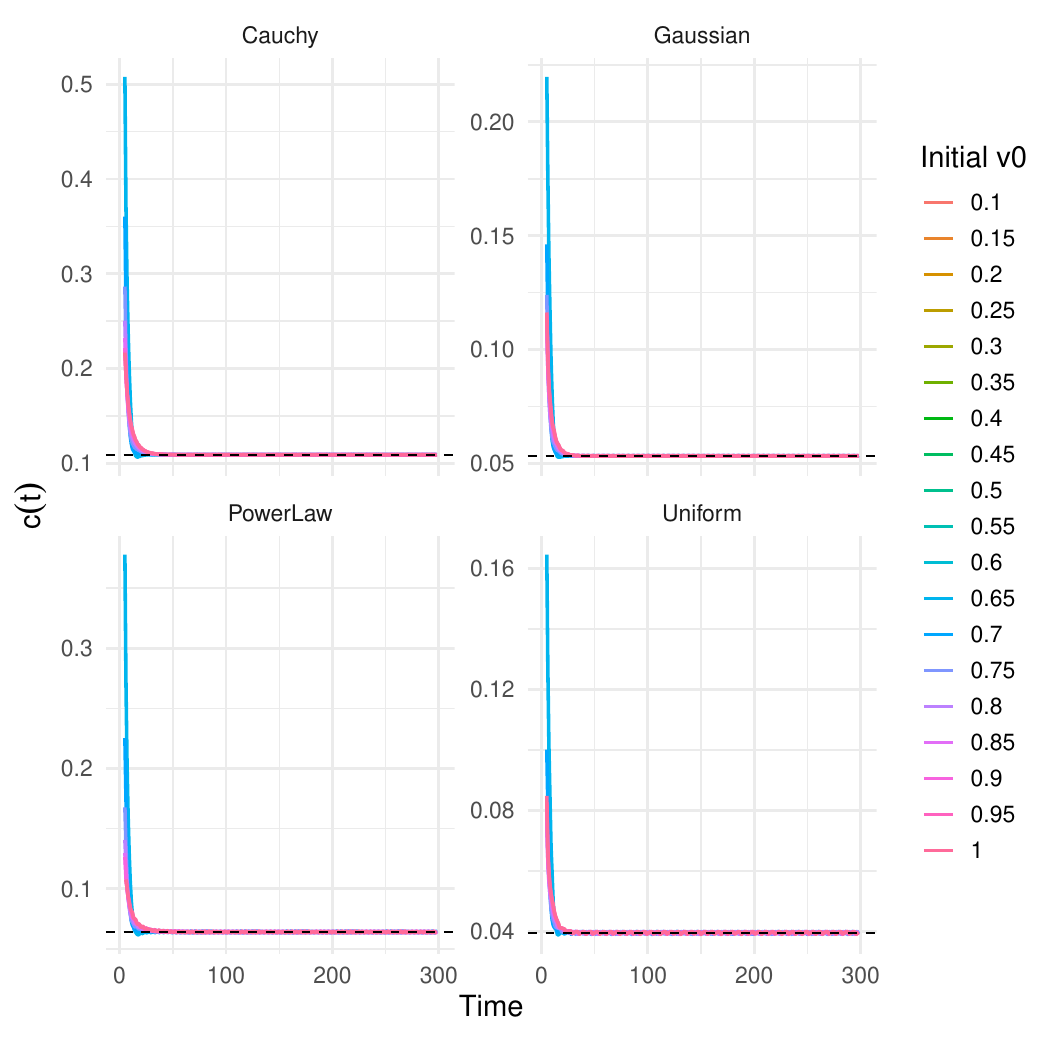} 
   \caption{Wave Speed $c(t)$ for Different Initial Releases $v_0$. The dashed line represents asymptotic wave speed for each kernel.  We have that $\lim_{t\to\infty }c(t)=0.1$ the Cauchy kernel, $\lim_{t\to\infty }c(t)=0.05$ for the Gaussian kernel, $\lim_{t\to\infty }c(t)=0.04$  for the Uniform kernel, $\lim_{t\to\infty }c(t)=0.06$ the Power Law kernel. These plot once again confirm that of the four kernels considered, the Cauchy dispersal kernel produces the greater wave speeds in general for different initial releases incidence $v_0$, followed respectively by the Power Law, Gaussian, and Uniform kernels.}
   \label{fig:example}
\end{figure}

\subsection*{Biological importance}
In the context of {\it Wolbachia} dynamics modeled by Lattice Difference Equations (LDEs), the \emph{wave speed} is a crucial quantity that quantifies how fast the infection front propagates through a mosquito population distributed in space. Understanding and predicting this speed has several important implications:

\begin{itemize}
    \item \textbf{Biological Insight:} The wave speed determines how quickly {\it Wolbachia} can spread spatially, affecting the feasibility of using it as a biological control agent against mosquito-borne diseases such as dengue, Zika, and malaria.

    \item \textbf{Invasion Success:} For a successful invasion, the wave speed must be positive, i.e., the infection must advance spatially over time. A negative or zero wave speed would indicate that {\it Wolbachia} fails to establish itself.

    \item \textbf{Threshold Phenomena:} In models with an \textit{Allee effect} or \textit{threshold dynamics}, the wave speed helps identify critical release strategies. Only above a certain infection threshold will a traveling wave with positive speed exist.

    \item \textbf{Optimization of Release Strategies:} Knowledge of the wave speed enables public health agencies to design optimal release programs by balancing release density, frequency, and spatial coverage to ensure the timely spread of {\it Wolbachia}.

    \item \textbf{Connection to Dispersal Kernels:} Different dispersal kernels (e.g., Gaussian, Cauchy, Uniform, Power-law) lead to different wave speeds. The wave speed thus reflects the interplay between local dynamics (via the {\it Wolbachia} growth function) and spatial dispersal mechanisms.

    \item \textbf{Predictive Modeling:} Accurately computed wave speeds allow researchers to forecast the temporal and spatial scales required for {\it Wolbachia} to dominate in heterogeneous or fragmented environments.

\end{itemize}

In summary, the wave speed is not only a mathematical descriptor but also a biologically relevant quantity that informs both the understanding and practical application of {\it Wolbachia}-based vector control programs. In general, heavy-tailed kernel such as the Cauchy and Laplacian produce higher wave speeds.

%%%%%%%%%%%%%%%%%%%%%%%%%%%%%%%%%%%%%%%%%%%%%%%%%%%%%%%%%%%%%%
\section{Distribution of outbreak size}\label{sec:infectionburder}
%%%%%%%%%%%%%%%%%%%%%%%%%%%%%%%%%%%%%%%%%%%%%%%%%%%%%%%%%%%%%%

\begin{figure}[H] 
  \resizebox{0.9\textwidth}{!}{\begin{minipage}{1.4\textwidth}
 \begin{tabular}{cc}\\
 {\bf $k=1$} & {\bf $k=3$}\\
  \includegraphics[width=4.5in]{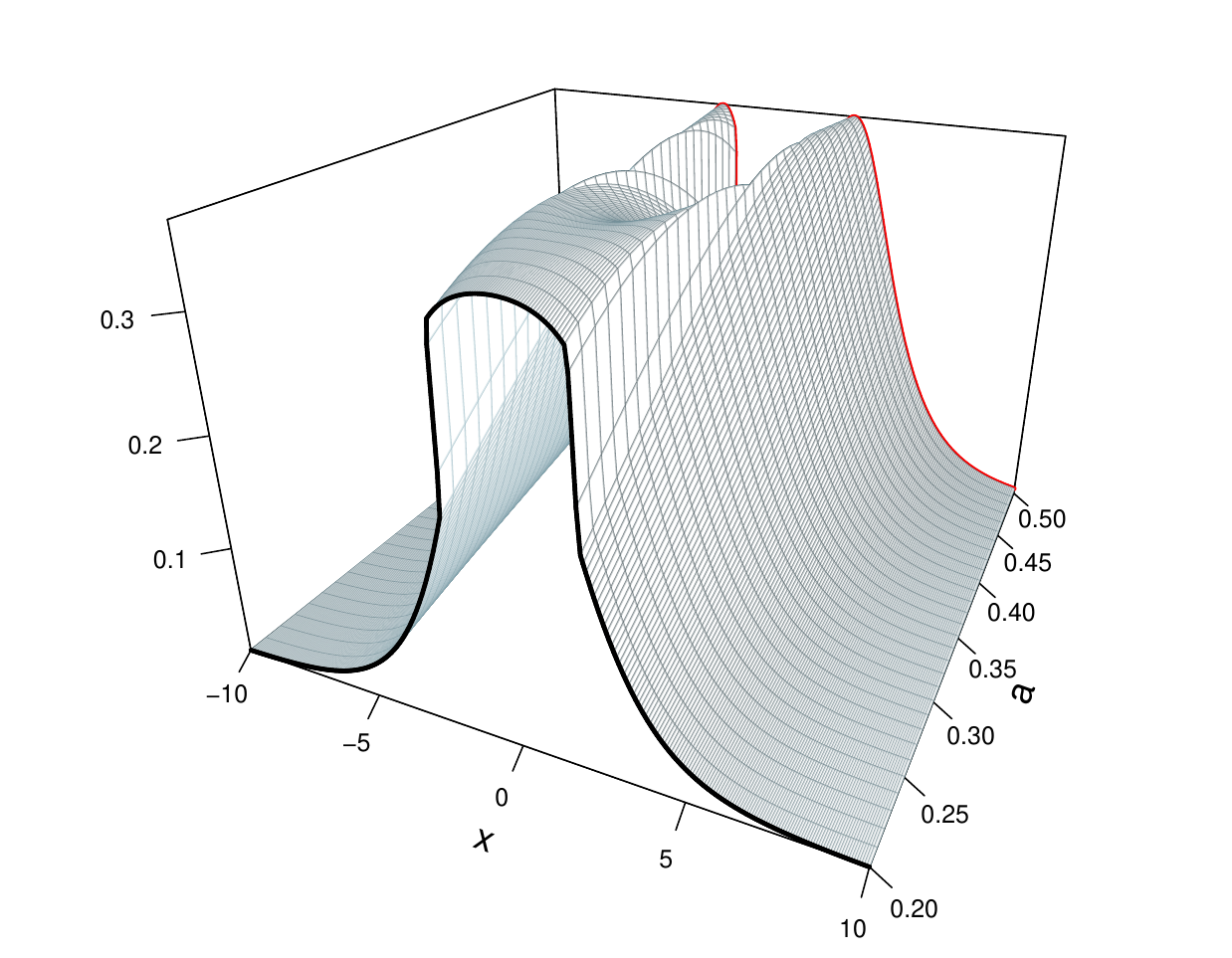}  &   \includegraphics[width=4.5in]{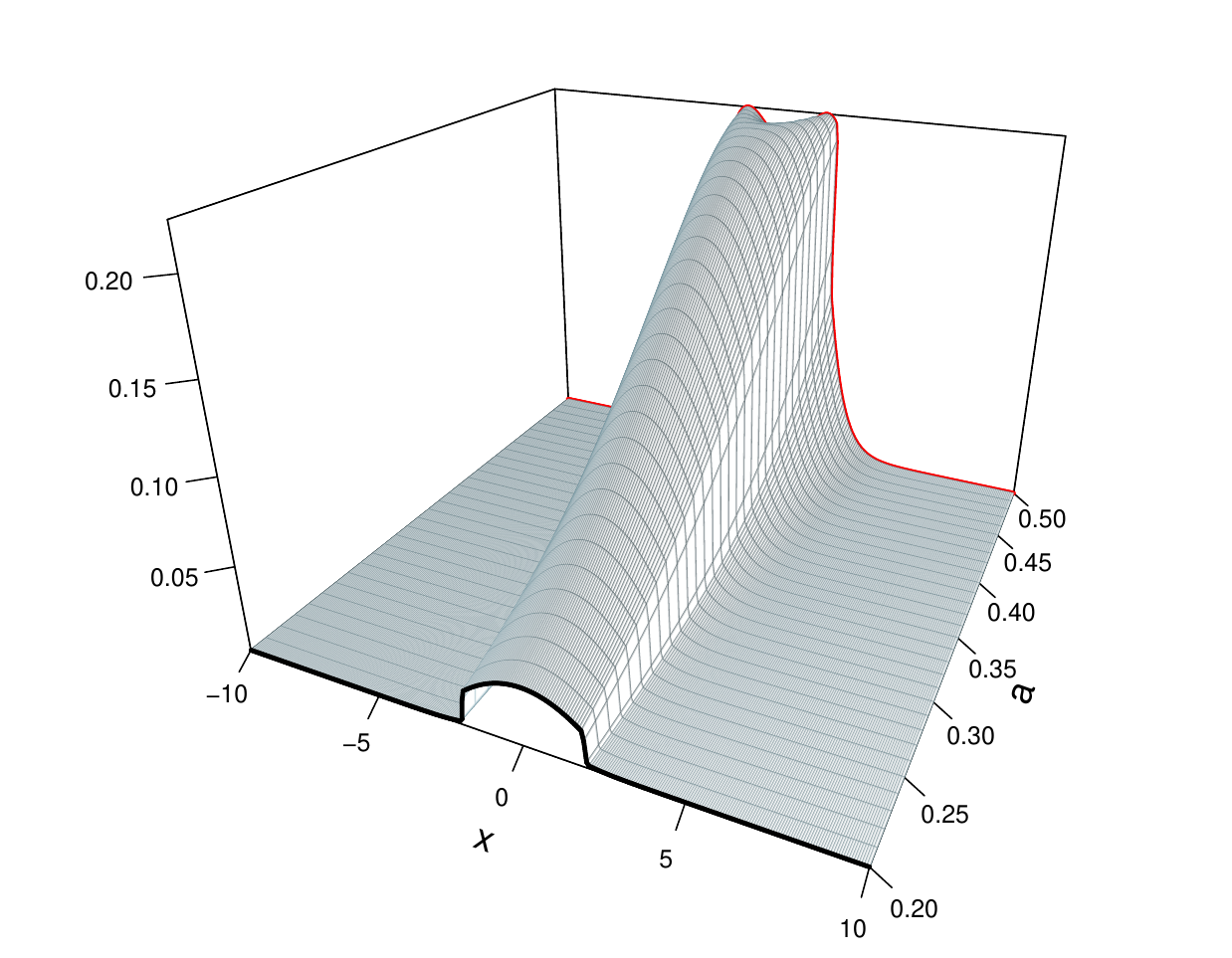} \\
   {\bf $k=10$} & {\bf $k=25$}\\
    \includegraphics[width=4.5in]{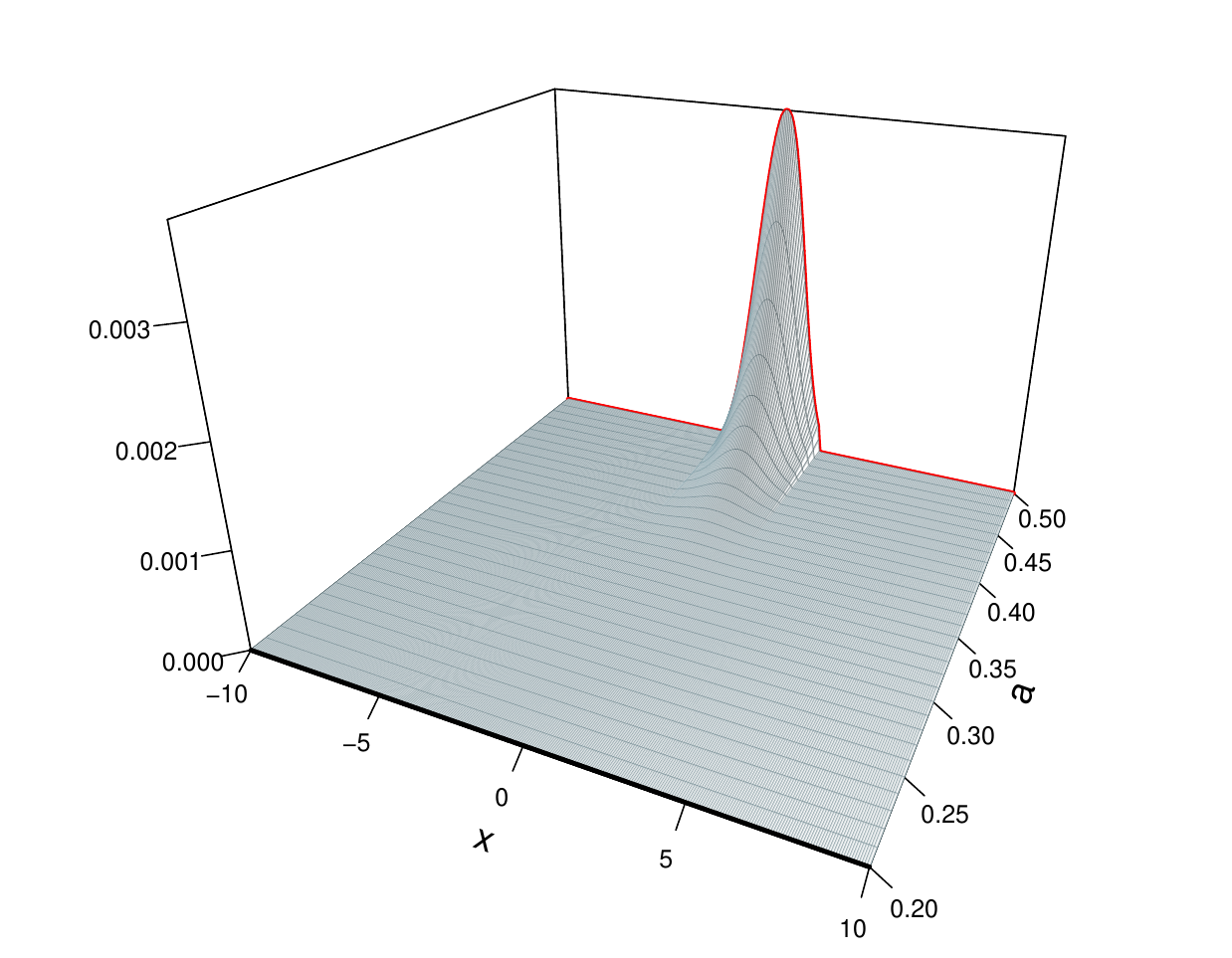}  &   \includegraphics[width=4.5in]{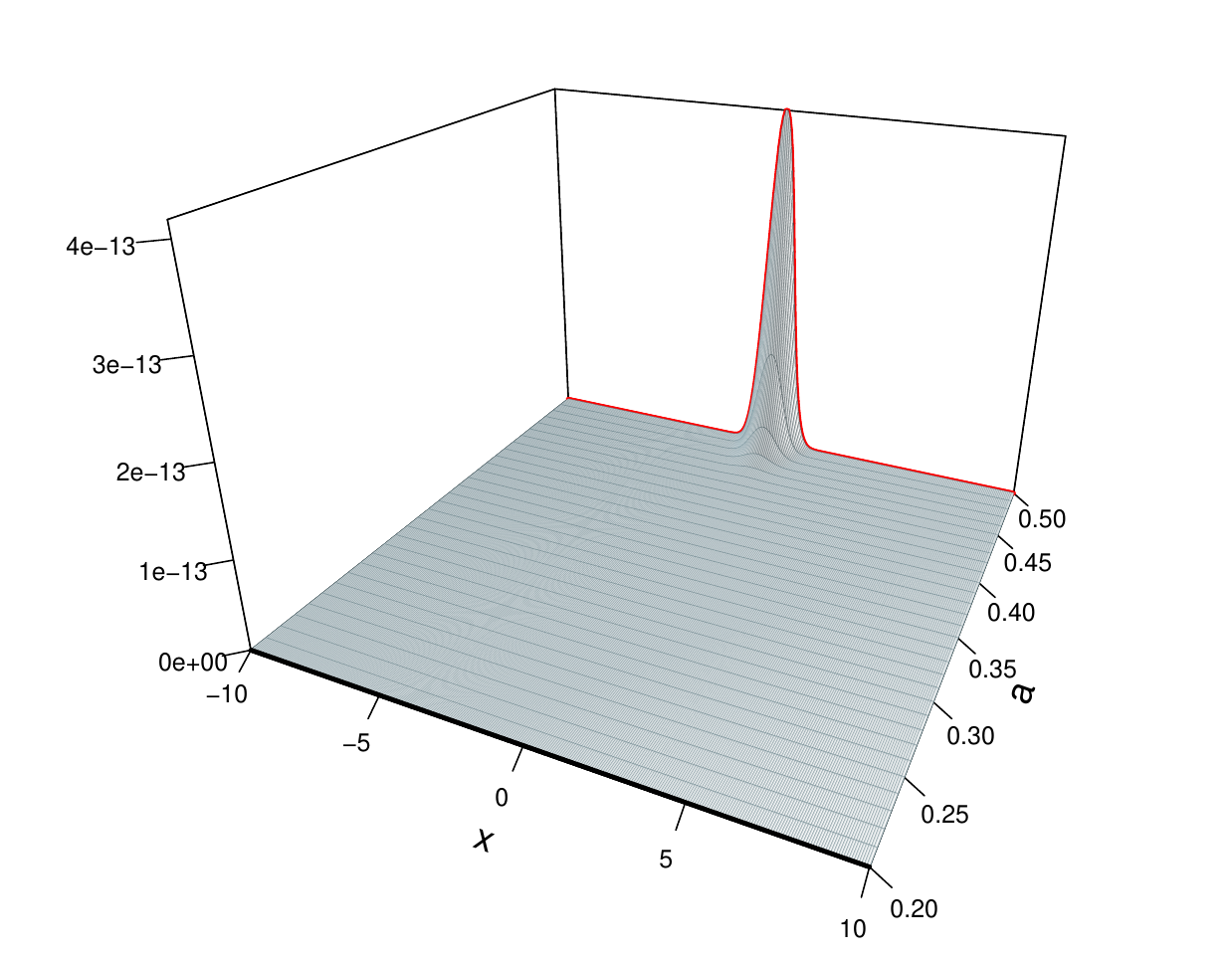} 
    \end{tabular}
    \end{minipage}}
    \caption{Distribution of outbreak size for size $k=1, 3, 10$ and $25$. The front curve (black) represents the initial distribution with $a=0.2$ and the back curve represents the final distribution (a=0.5) using a Gaussian kernel.}
          \label{fig:IB-Gaussian}
 \end{figure}

There is an interesting observation when the {\it Wolbachia} growth function is used. We recall that $v(t,x)$ represents the infection frequency at location $x$ and time $t$. In the sequel, a specific location will be denoted  as $x_i$ to stress the importance of the local nature of the results that will follow. From the derivation of $v(t,x_i)$ in \cite{Yu2019}, it follows that it represents the ratio of infected female or males relatively to the whole population at time $t$ and location $x_i$. As such, it can be understood as the probability of an infection at time $t$ and location $x_i$. This allows us to define the following random variables: given a location $x_i\in \R$ and time $t\in \Z^+$, 
\[
Y_{it}=\begin{cases}
1 & \textrm{if there is an infected individual at location $x_i$ and time $t$}\\
0 & \textrm{otherwise}
\end{cases}\;.
\]
Therefore, we have that 
\[
\begin{cases}
\Pp(Y_{it}=1)=v(t,x_i)\\
 \Pp(Y_{it}=0)=1-v(t,x_i)
\end{cases}\;.
\]
Clearly, if $v(t,x_i)=0$, then there is no individual infected at location $x_i$ at time $t$, therefore, $\Pp(Y_{it}=1)=0$ and $\Pp(Y_{it}=0)=1$.
Likewise, if $v(t,x_i)=1$, then there is an individual infected at location $x_i$ at time $t$, therefore, $\Pp(Y_{it}=1)=1$ and $\Pp(Y_{it}=0)=0$.
We define, for all $x_i\in \R$
\[
N_i=\min\set{t\in \Z^+: v(t+1,x_i)=v(t,x_i)}, \quad \mbox{$p_{it}:=v(t,x_i)$ for  $0\leq t\leq N_i-1$}\;.
\]
Let $\ds Y_i|N_i=\sum_{t=0}^{N_i-1}Y_{it}$ be the number of infected individuals at location $x_i$ from time $t=0$ to $t=N_i-1$. We will refer to $Y_i$ as the outbreak size. 
We have the following result:
\begin{thm}\label{thm:disOutbreak}
Given $N_i$ as above, $Y_i|N_i$ has a Poisson-Binomial distribution, that is, 
\[
Y_i|N_i=\sum_{t=0}^{N_i-1}Y_{it}\sim PB(p_{i0}, p_{i1},\cdots, p_{i,N_i-1})\;.
\]
More over, there exists $0<q_i<1$ such that for any $k\in \Z^+$, 
\begin{equation}\label{eqn:distOfOutbreal}
\Pp(Y_i=k)=\sum_{m=0}^{\infty}\Pp(Y_i=k|N_i=m)(1-q_i)^mq_i\;,
\end{equation}
where
 \begin{equation}\label{eqn:CondDistOutbreak}
 \Pp(Y_i=k|N_i=m)=\sum_{U\in \Lambda_k}\prod_{j\in U}p_{ij}\prod_{i\in U^c}(1-p_{ij}) \;,
 \end{equation}
and 
 \[\Lambda_k=\set{S_k:  S_k\subseteq \set{1,2,\cdots, m}, ~|S_k|=k}, \quad \mbox{ $|A|$ is the cardinality of the set $A$}\;.\]
\end{thm}
The proof can be found in Appendix C.
 The  probabilities $\Pp(Y_i=k)$  for $k\in \Z^+$ define the probability mass function of outbreak size distribution and thus  will give spatial insight into how likely specific outbreak sizes are under the {\it Wolbachia} invasion scenario.  In general, the expression in \eqref{eqn:CondDistOutbreak} is very complex to obtain because given $k$, $\Lambda_k$ contains $k \choose{m}$ values. Ultimately, the sum above is over $2^m-1$ values, which grows fast for larger values of $m$. For simplicity of illustration below, we will will fix $N_i=400$ and assume that the $p_{it}$'s are small, so that 
 $Y_i\sim Poisson(\lambda_i)$ with $\lambda_i=\sum_{t=0}^{N_i-1}p_{it}$. In this case
  \begin{equation*}
 \Pp(Y_i=k)=e^{-\lambda_i} \frac{\lambda_i^k}{k!}\;.
 \end{equation*}
 The starting infection profile $v_0(x)=v(x,0)$ plays an  important role in understanding the distribution of outbreak sizes. For simplicity, we will choose a pulse profile function  $v_0(x)=a\cdot 1_{\set{|x|<L}}(x)$, for a fixed value  $L=2$ and different equidistant values of $a$ between 0.2 and 0.5, where $1_A(x)$ is the characteristic function of the set $A$. As the amplitude of the function $v_0(x)$, $a$ represents the initial infection rate over the domain $[-L,L]$. In the figures below, we plot the likelihood of outbreak size $Pr(Y_i=k)$ for $k=1, 3, 10, 25$, given that $s_f= 0.3, s_h=0.7,\delta =0.5$ for respectively a Gaussian kernel (Figure \ref{fig:IB-Gaussian}) and Laplace kernel (Figure \ref{fig:IB-Laplace}). For $k=1$, Figure \ref{fig:IB-Gaussian} shows that  the likelihood of finding an infected individual is higher at the release point $x=0$ for $a=0.2$. As $a$ increases from 0.2 to 0.5, that likelihood (($\Pp(Y_i=1)> 0.3$)) is  the same for now two points somewhat far ($x_i\approx -5$ and $x_i \approx 5$) from the initial release point  but decreases significantly at the initial release point. For $k=3$, the likelihood  of finding three infected individual is low ($\Pp(Y_i=3)< 0.05$) around the release point but increases with $a$ ($\Pp(Y_i=3)> 0.2$). For $k=10$ and $k=25$,  Figure \ref{fig:IB-Gaussian} shows that the likelihood of finding an infected individual anywhere in space grows with $a$ from impossible ($\Pp(Y_i=10,25)\approx 0$) to very unlikely ($\Pp(Y_i=10,25)<0.003$). Using a Cauchy's  kernel and  Figure  \ref{fig:IB-Laplace}, we arrive at similar conclusions as above. 
 \begin{remark}
 We note that when $N_i\leq t\leq M$ for some $M\in Z^+$,  we have that 
 \[Y_i|N_i\sim Poisson(\lambda), \]
 where $\ds \lambda=\sum_{t=N_i}^Mp_{i*}=(M+1-N_i)p_{i*}$ with $p_{i*}=v^*=\frac{s_f}{s_h}$ is the unstable fixed point.\\
 Therefore, $Y_i|N_i$ does not converge almost surely to any discrete distribution if $M\to \infty$.   
 \end{remark}
 
  \begin{figure}[H] 
  \resizebox{1\textwidth}{!}{\begin{minipage}{1.4\textwidth}
 \begin{tabular}{cc}\\
 {\bf $k=1$} & {\bf $k=3$}\\
  \includegraphics[width=4.5in]{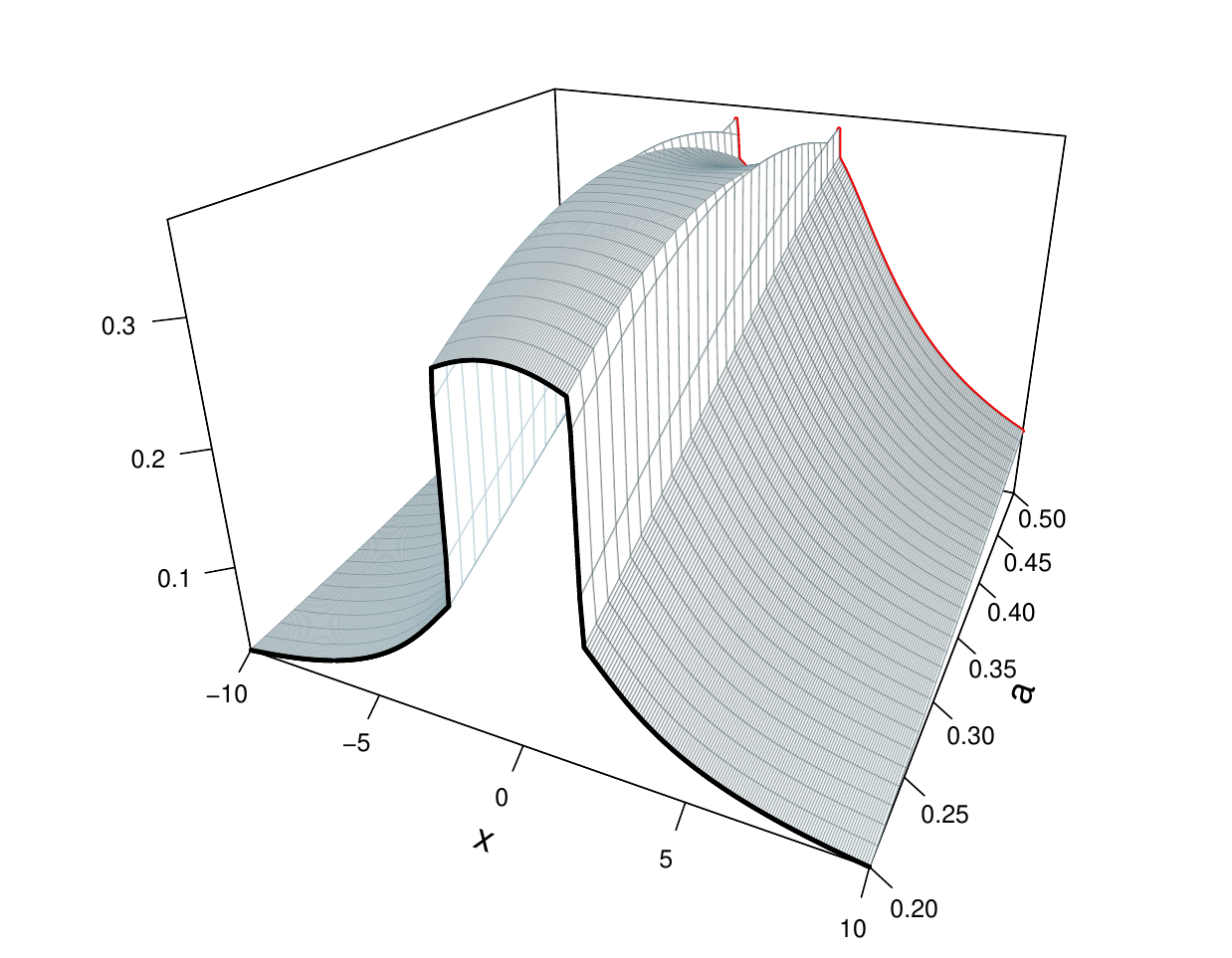}  &   \includegraphics[width=4.5in]{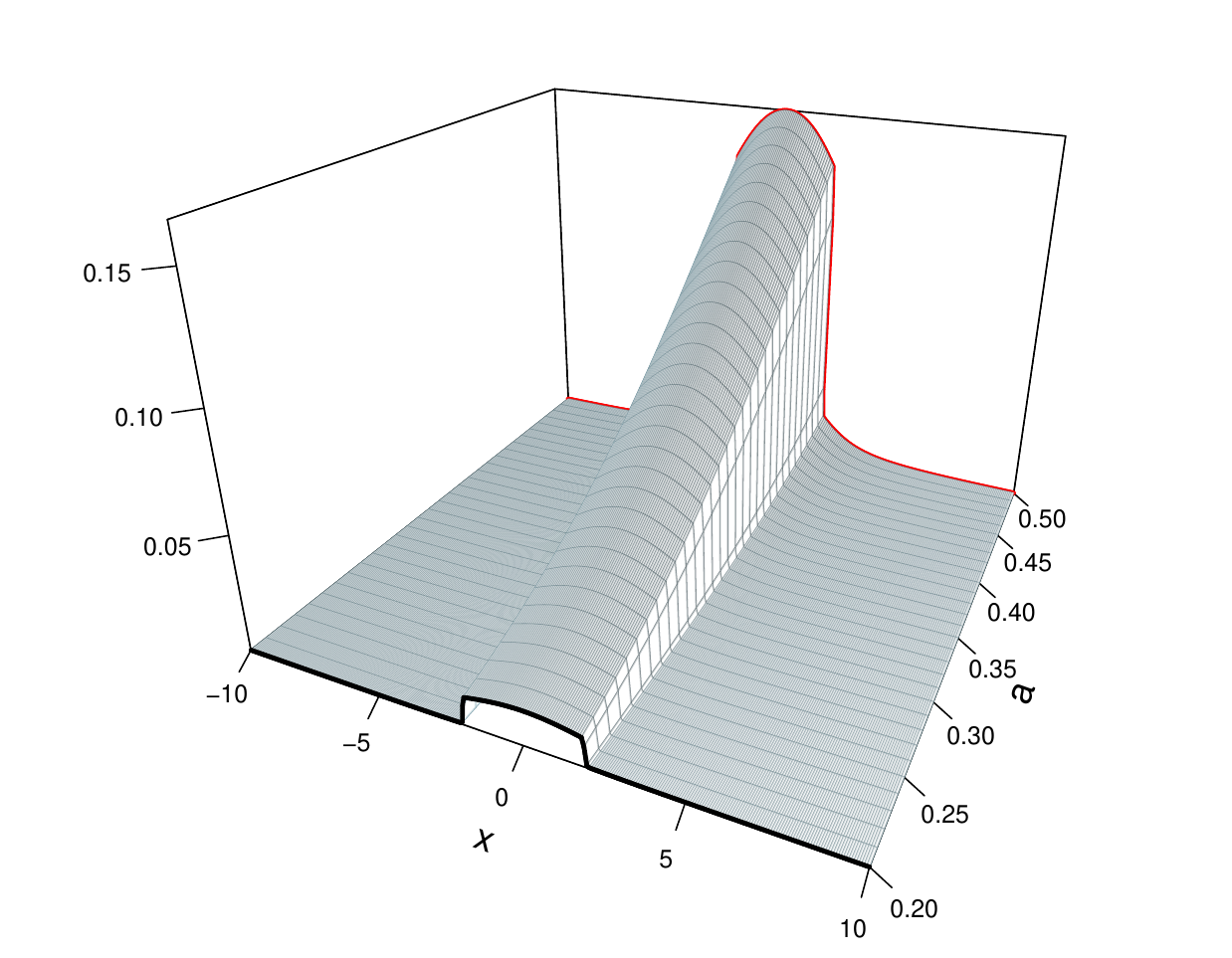} \\
   {\bf $k=10$} & {\bf $k=25$}\\
    \includegraphics[width=4.5in]{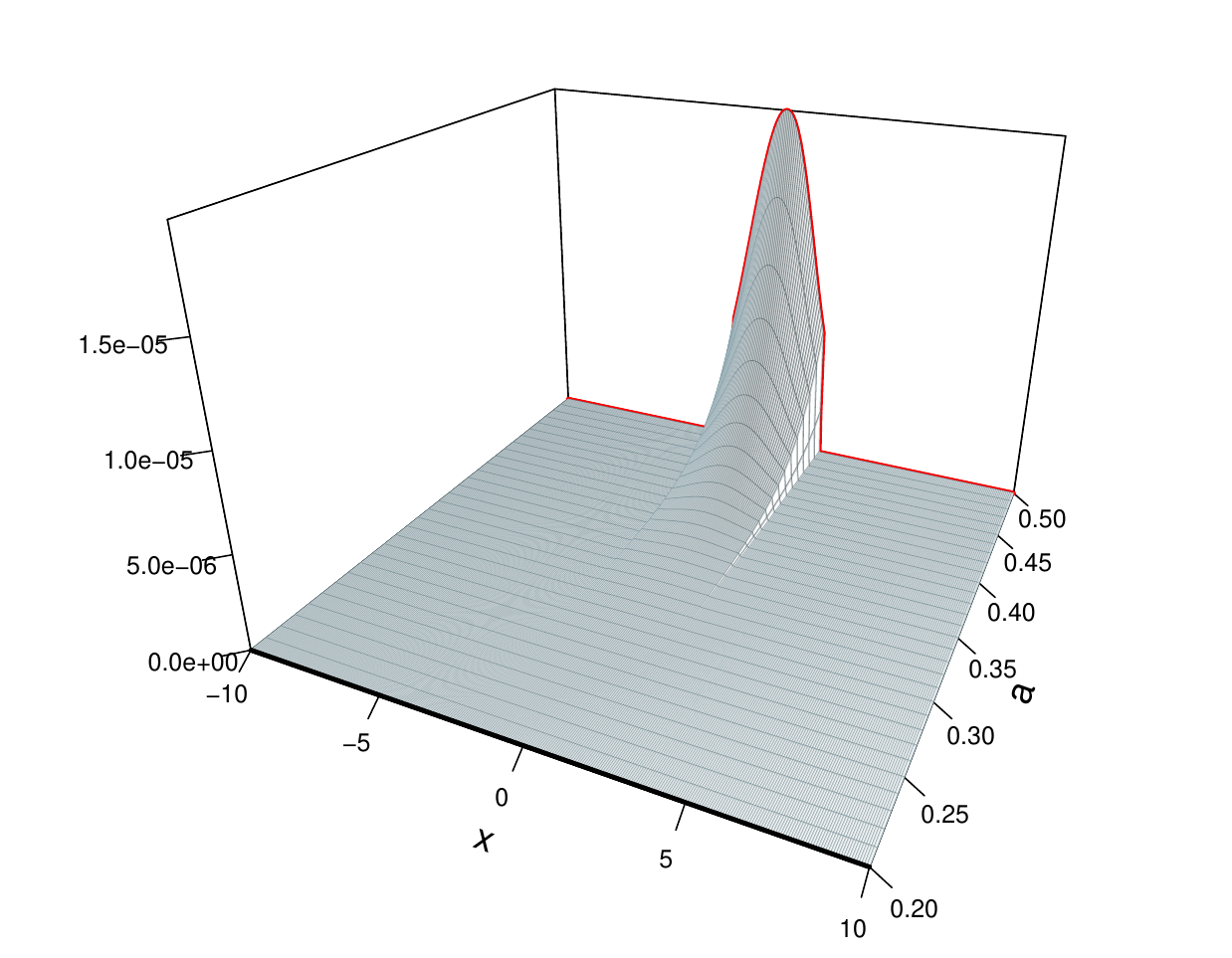}  &   \includegraphics[width=4.5in]{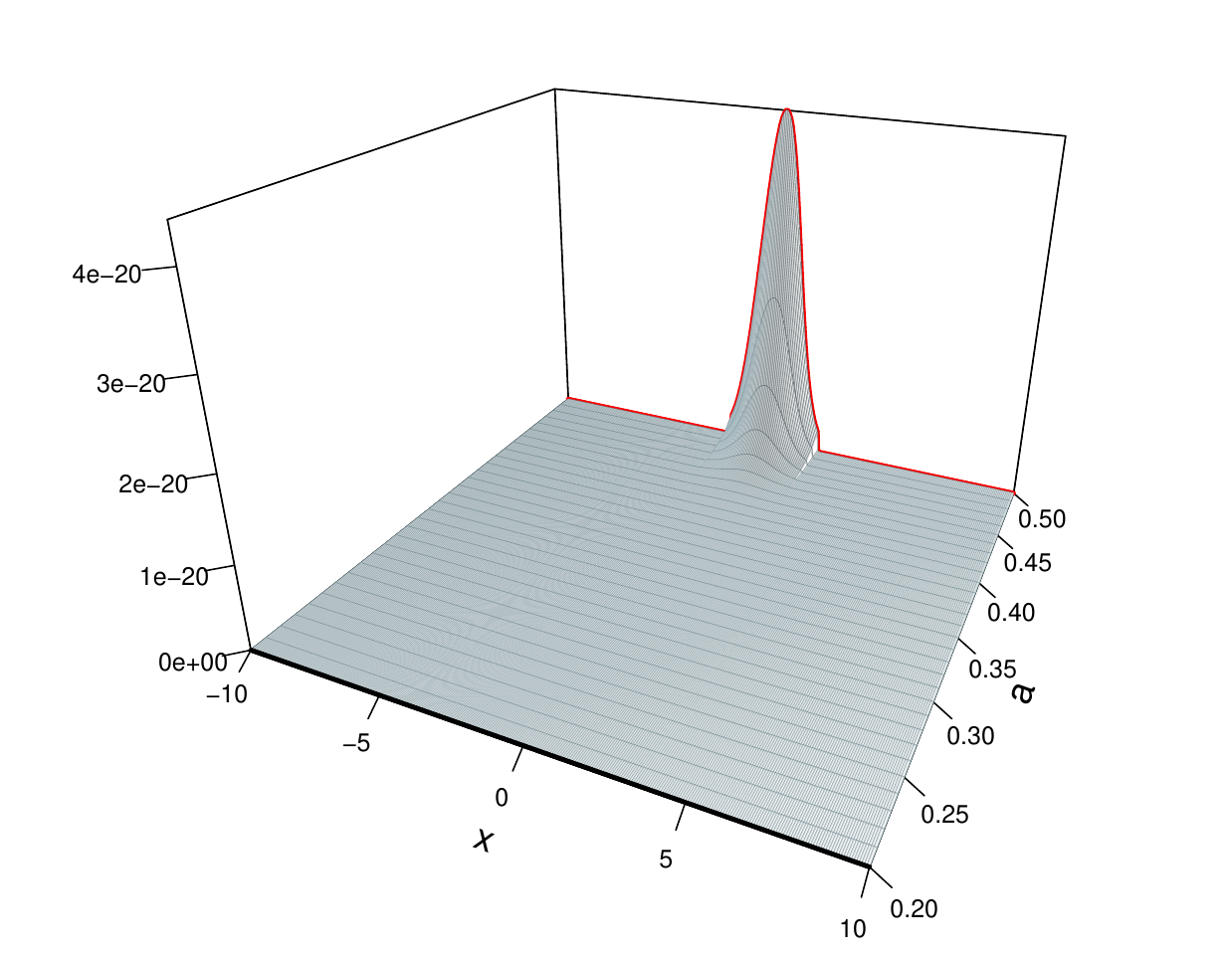} 
    \end{tabular}
    \end{minipage}}
       \caption{Distribution of outbreak size for size $k=1, 3, 10$ and $25$. The front curve (black) represents the initial distribution with $a=0.2$ and the back curve represents the final distribution ($a=0.5$) using a Laplacian kernel.}
     \label{fig:IB-Laplace}
 \end{figure}
 \subsection{Bi-modality as cost minimization criterion}
 What the previous simulations reveal is that at low outbreak sizes ($k$ small), the spatial distribution of outbreak sizes ($\Pp(Y_i=k))$)  has two regimes: a unimodal regime which corresponds to low values of the amplitude $a$ and a bimodal regime which  corresponds to high values of $a$. For low values of $a$, regions around the initial release point $x=0$ inside the release interval $[-L,L]$ are very likely to be infected whereas regions outside this interval are unlikely to be infected overtime. However, larger values  of  the amplitude $a$ create situations where spatially, the initial release  point and interval are less likely to be infected than certain regions outside of this interval. That change of modality occurs for a critical value, say $a^*$. This  critical value can be obtained by changing values of $a$ and tracking the mode of the distribution until it switches from unimodal to bimodal.  Dynamically, this may well be a point of bifurcation in the underlying dynamical system. We suspect  that this may be due to the presence of the Allee threshold in the growth function $f$.
  More importantly, this critical value can be used to compare different initial release profiles $v_0(x)$. To clarify this point further, we will consider the following three initial  release profiles, whose expressions are given in Table \ref{tab:InitialProfiles}, illustrated in Figure \ref{fig:InitialProfiles} below. We observe that the pulse profile provides the lowest critical value ($a^*=0.250$) followed by the Quadratic ($a^*=0.290$) and then the Triangular ($a^*=0.330$) profiles. Moreover, $a^*$ represents the minimum value of $a$ for which the infection spread to other parts of the domain. This makes it a useful early warning signal for impending regime shifts. 
  
   \begin{table}[h]
\resizebox{1\textwidth}{!}{\begin{minipage}{1\textwidth}
\centering \begin{tabular}{|l|l|l|}
\hline
{\bf Release Profile function} &Expression & Cost \\ \hline 
{\bf Pulse} &$v_0(x)=a\cdot 1_{_{|x|\leq L}}(x)$ & $2aL$ \\  \hline 
{\bf Triangular}&$v_0(x)=\ds a\left(1-\frac{|x|}{L}\right) \cdot 1_{_{|x|\leq L}}(x)$ & $aL$\\ \hline
{\bf Quadratic} &$v_0(x)=\ds a\left(1-\frac{x^2}{L^2}\right)\cdot 1_{_{|x|\leq L}}(x)$ & $\ds \frac{4aL}{3}$\\ \hline
\end{tabular}
\end{minipage}}
\caption{Selected initial release functions and their costs.}
\label{tab:InitialProfiles}
\end{table}
 \begin{figure}[H] 
    \centering
    \includegraphics[width=5in]{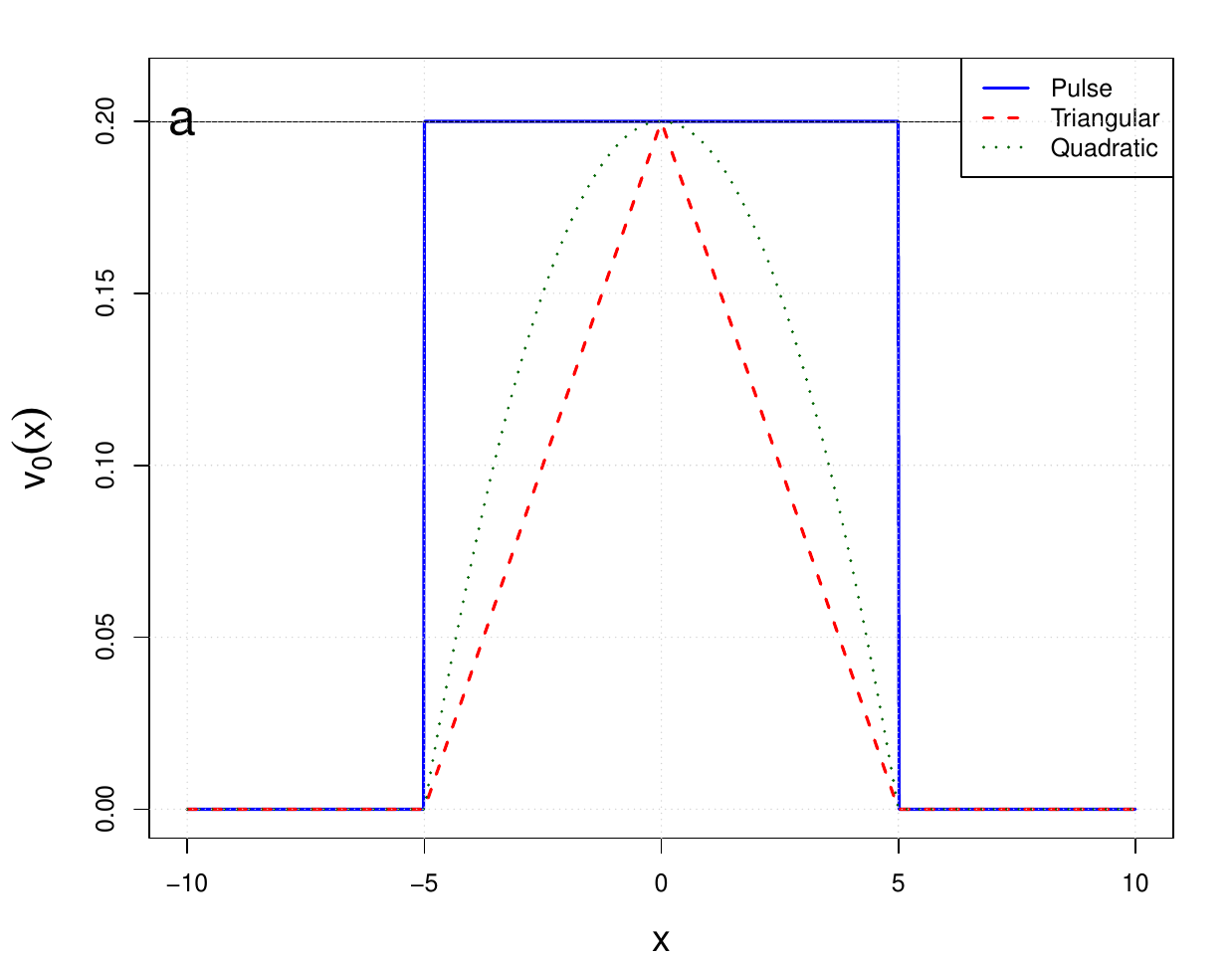} 
    \caption{An illustration of the  initial release functions}
    \label{fig:InitialProfiles}
 \end{figure}
 
 In the figures below, obtained for the same parameters as above and using a Lalacian  Kernel, we compare the probability of outbreak size $k=1$ for initial release profiles above. We find the for each profiles, the minimal value $a^*$ of $a$ or critical value  for the switch between a unimodal regime and a bimodal regime.
 \begin{figure}[H] 
  \resizebox{1\textwidth}{!}{\begin{minipage}{1.4\textwidth}
 \begin{tabular}{cc}\\
 {\bf Pulse} & {\bf Triangular}\\
  \includegraphics[width=4.5in]{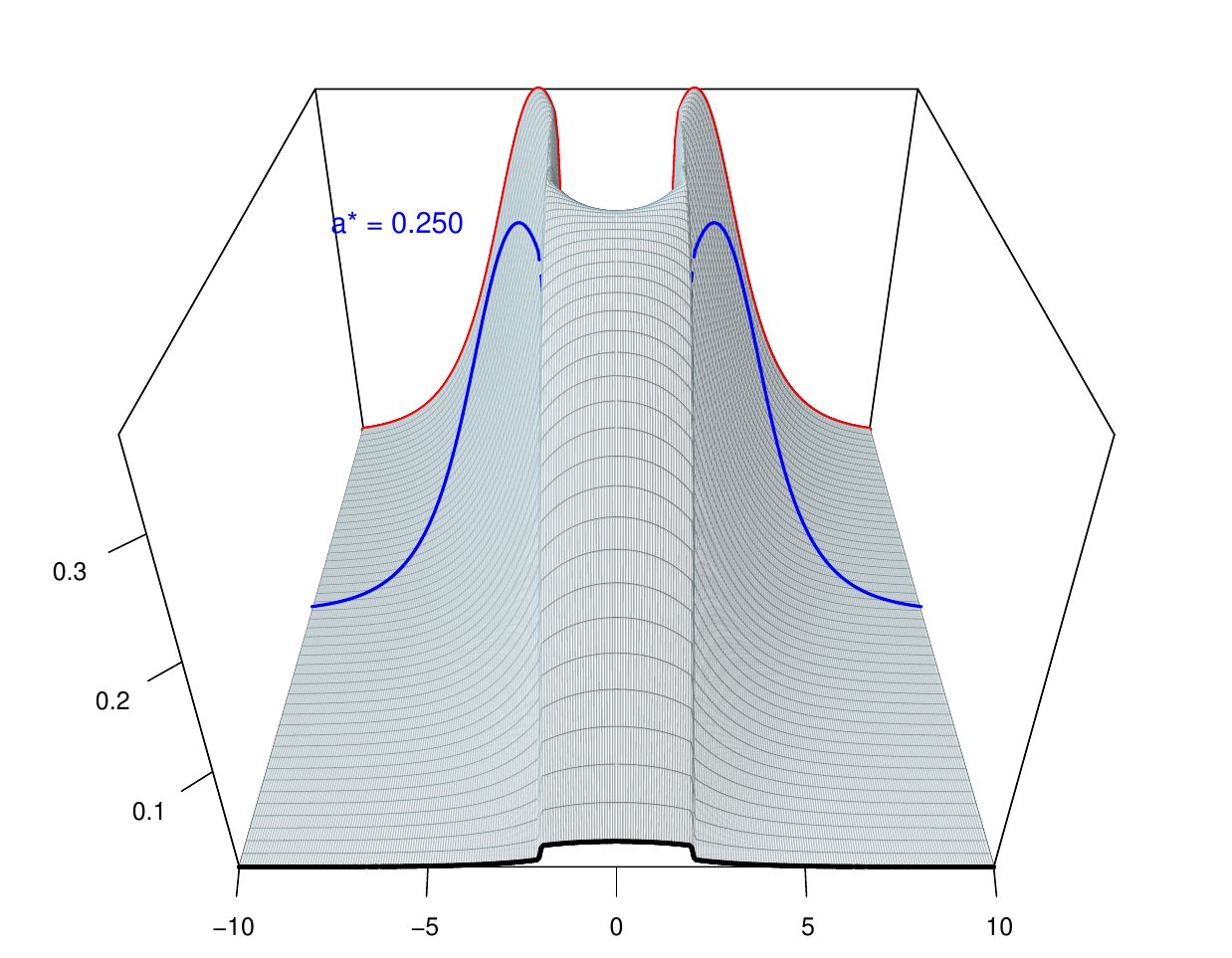}  &   \includegraphics[width=4.5in]{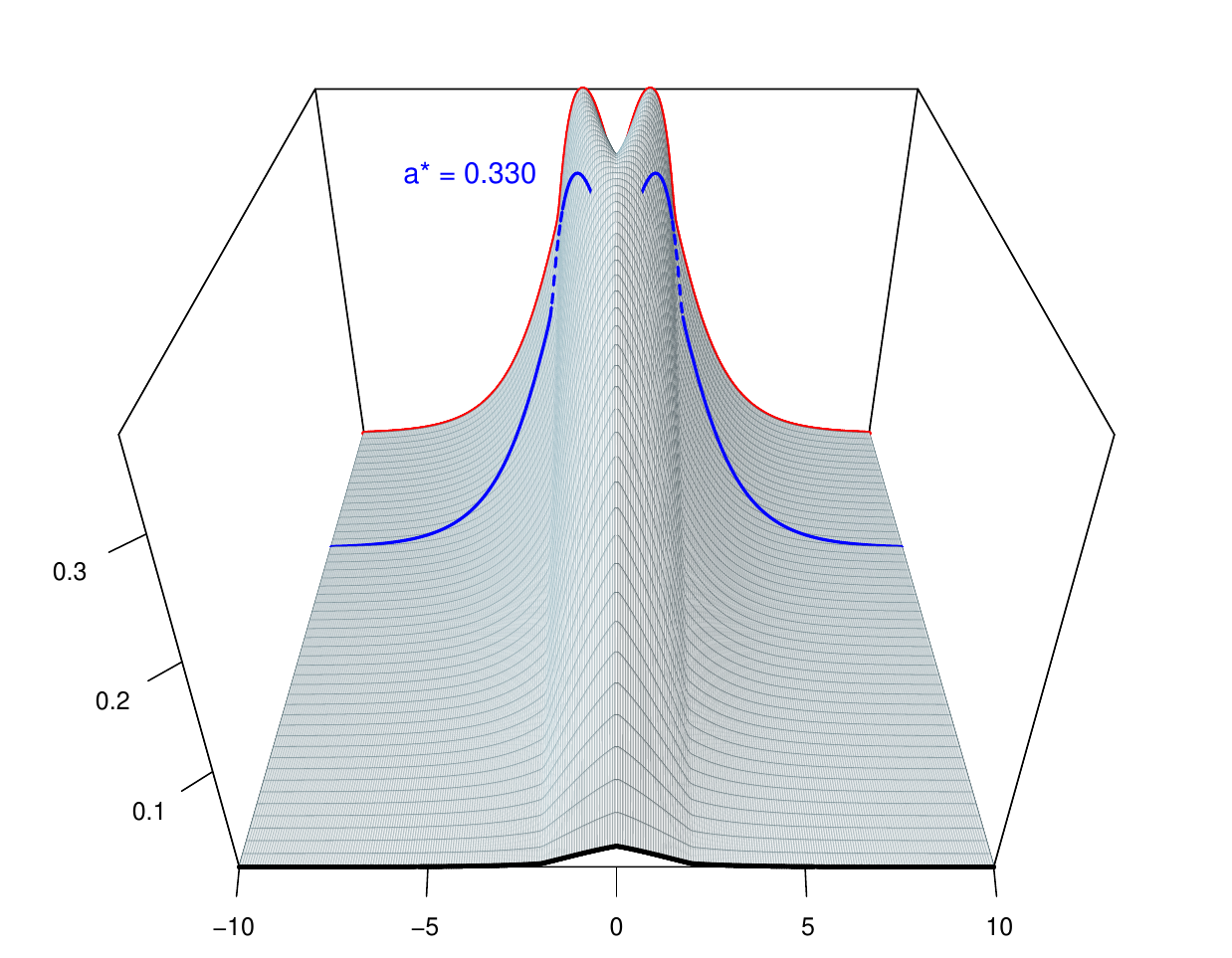} \\
   {\bf Quadratic} & \\
    \includegraphics[width=4.5in]{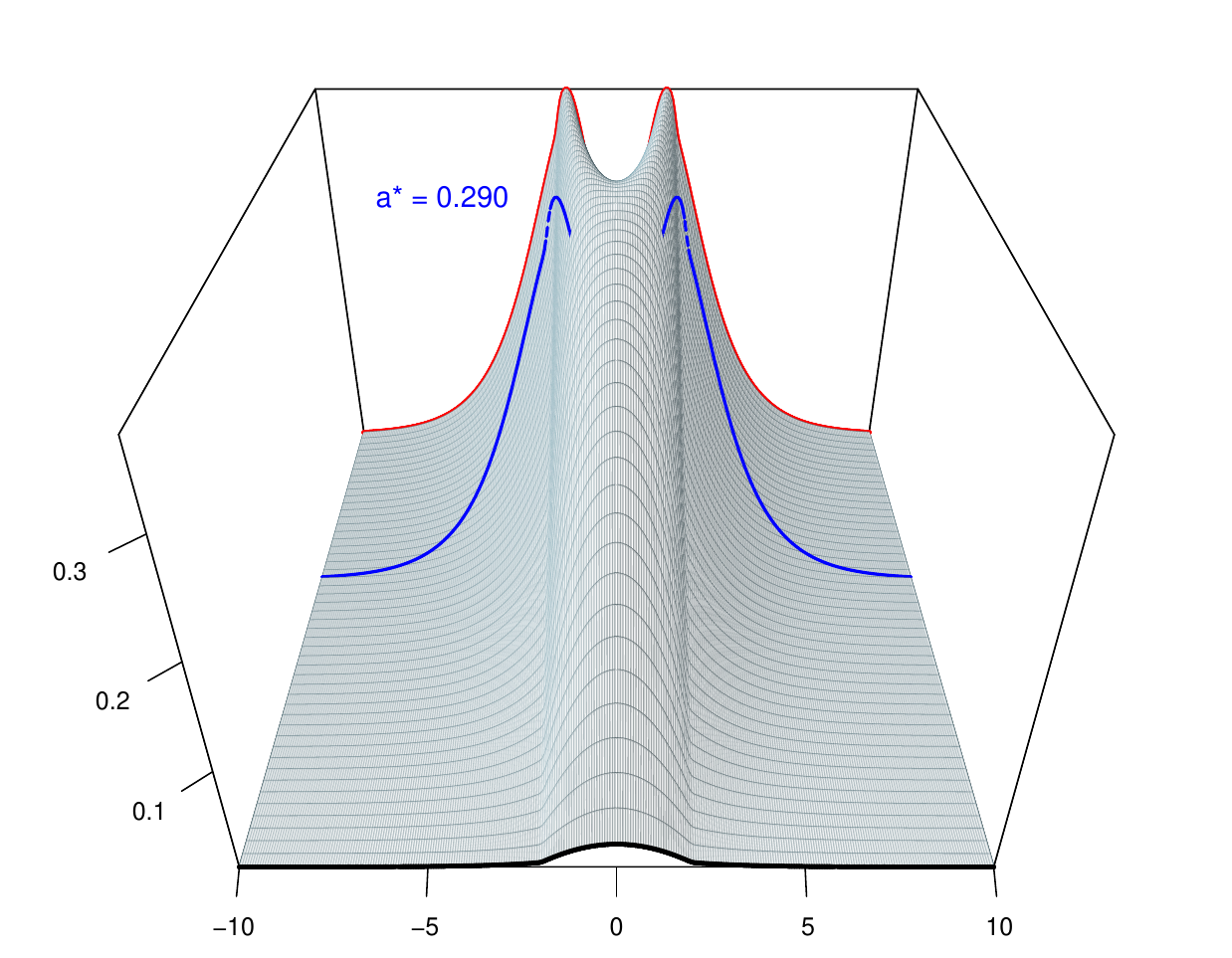}  &  
    \end{tabular}
    \end{minipage}}
       \caption{Distribution of outbreak size for size $k=1$ obtained for different initial release profiles. The blue curve represents the distribution of outbreak size $k=1$, when $a$ has reached the critical values at which the distribution becomes bimodal. }
     \label{fig:IB-Laplace}
 \end{figure}
   
  For release strategies, finding this minimal  value of $a$  in conjunction with $L$ could help minimize the cost of release (see \cite{Turelli2017Wolbachia, Zhou2008Stability, Chan2020Control}), that is, the total biological or logistical effort required to initiate an infection invasion by releasing {\it Wolbachia}-infected individuals (e.g., mosquitoes) into a population. It is  given as $Cost (a,L)=\int v_0(x)dx$. This quantity reflects the proportion of mosquitoes infected at time $t=0$, accounting for both amplitude and spatial extent. Classical minimization (Asymptotic Constraints Minimization or ACM) of cost techniques  (see \cite{Chan2020Control, Muratov2012Bistable, Benguria1996FrontSpeed}) use asymptotic constraints, that is, 
  \[
(a^*,L^*)=\mbox{Argmin}\set{Cost(a,L) \quad s.t. ~~v(0,x)=v_0(a,L)~\quad \mbox{and $~\ds \lim_{t\to \infty} v(t,x)=1$}}\;.
  \]
  Using bi-modality, we propose a minimization criterion, that we call Modality Constraint Minimization (MCM),  based on the first switch to a bimodal regime. Noting that the mode  of the likelihood of outbreak release $\Pp(Y_i=k)$ depends on the amplitude $a$ and the width  $L$ of the initial  release interval, we formality state this binary minimization as 
    \[
 (a^*,L^*)= \mbox{Argmin}\set{Cost (a,L) \quad s.t. ~~ \mbox{Mode}(a, L)=2}\;.
  \]
  Let us offer a few differentiators of the two methods. 
  \subsubsection*{Advantages of MCM}
  1. ACM assumes perfect deterministic knowledge of the system whereas MCM handles uncertainty (heterogeneity of landscapes or discrete values of $k$).\\
  %2. If the system is stochastic or parameter-sensitive, MCM gives a robust probabilistic threshold rather than an idealized one.\\
  2. ACM requires full convergences $\ds \lim_{t \to \infty}v(t,x)=1$ which is expensive to verify and maybe unnecessary if high probability of invasion would suffice. MCM replaces a hard convergence constraint with a statistical phase transition, that is, the emergence of bi-modality signals a non-negligible chance of invasion.\\
  3. MCM inherently reflects distributional shifts, thus captures critical transitions. This makes MCM interpretable via early warnings which ACM cannot. 
    \subsubsection*{Limitations of MCM}
1. MCM only signals mixed outcome probability, not invasion certainty. For example, the system can be in bimodal regime with only a $30\%$ change of successful invasion.\\
2. Bi-modality detection can be noisy or resolution-dependent. It requires a good tuning of parameters and needs a well-defined mode separation criteria.\\
3. MCM needs many ensemble simulations to estimate outcome distributions reliably. This may become computationally expensive in the case of high-resolution spatial model.

In all, the advantages of MCM do not make it better than ACM nor do its limitations make it worse. MCM does  offer a flexible approach and potentially more informative approach under uncertainty and could be used complementary to ACM. 

  \subsubsection*{Numerical comparison}
To have a complete picture that takes into account the kernel type, the initial profiles per outbreak size, while keeping other parameters as  $\delta=0.2, s_f=0.2, s_h=0.9, \beta=0.9, N_x=400, N_t=200 $, we propose  the  comparison Table \ref{tab:table4} below.  We note that the optimal values of $a$ and $L$ from the ACM are obtained using the bisection method (see \cite{burden2010numerical}), with threshold $\beta=0.9$.  The MCM turns out to be  very flexible (different optimal values based on outbreak size $k$). Regardless of the kernel, the Triangular profile is the most cost effective for both methods. We note that ACM is much more sensitive to the tuning of parameters, especially $s_f, s_h$, and $\delta$. If for instance we choose $s_f=0.3, S_h=0.7$ and $\delta=0.5$, then the optimal value of $a$ is $a^*=1$ for all kernels (for sake of brevity, we did not include it). This is not only biologically implausible, but also strategically costly. For the same values, the MCM produces somewhat  high values of $a^*$, but still significantly lower than 1.

\begin{table}[H]
\centering
\begin{tabular}{|c|c|c|c|c|c|c|}
\hline
\textbf{Kernel} & \textbf{Profile} & \textbf{Outbreak Size} & MCM \(a^*\) & ACM \(a^*\) & MCM Cost & ACM Cost \\
\hline
\multirow{12}{*}{Laplacian} 
& \multirow{4}{*}{Pulse} & $k=1$ & 0.200 & \multirow{4}{*}{0.395} & 0.200 & \multirow{4}{*}{0.395} \\
 &  & $k=2$ & 0.290 &  & 0.290 &  \\
 &  & $k=3$ & 0.340 & & 0.340 &  \\
 &  & $k=4$ & 0.360 & & 0.360 &  \\
%\cline{2-7} 
\hline
& \multirow{4}{*}{ Quadratic} & $k=1$ & 0.220 & \multirow{4}{*}{0.518} & 0.147 & \multirow{4}{*}{0.345} \\
 &  & $k=2$ & 0.320 &  & 0.213 &  \\
 &  & $k=3$ & 0.380 & & 0.253 &  \\
 &  &$k=4$ & 0.420 &  & 0.280 &  \\
%\cline{2-7}
\hline
& \multirow{4}{*}{Triangular} & $k=1$ & 0.260 & \multirow{4}{*}{0.720} & 0.130 & \multirow{4}{*}{0.360} \\
 &  & $k=2$ & 0.380 & & 0.190 &  \\
 &  & $k=3$ & 0.450 &  & 0.225 &  \\
 &  & $k=4$ & 0.490 &  & 0.245 &  \\
\hline \hline
\multirow{12}{*}{Gaussian} &  \multirow{4}{*}{Pulse} & $k=1$ & 0.200 & \multirow{4}{*}{0.390} & 0.200 & \multirow{4}{*}{0.390} \\
 & & $k=2$ & 0.290 & & 0.290 &  \\
 & & $k=3$ & 0.330 &  & 0.330 &  \\
 & & $k=4$ & 0.360 &  & 0.360 &  \\
%\cline{2-7}
\hline
&  \multirow{4}{*}{Quadratic} & $k=1$ & 0.220 & \multirow{4}{*}{0.493} & 0.147 & \multirow{4}{*}{0.329} \\
 & & $k=2$ & 0.330 &  & 0.220 & \\
 & & $k=3$ & 0.380 &  & 0.253 & \\
 & &$k=4$ & 0.420 & & 0.280 &  \\
%\cline{2-7}
\hline
&  \multirow{4}{*}{Triangular} & $k=1$ & 0.260 & \multirow{4}{*}{0.640} & 0.130 & \multirow{4}{*}{0.320} \\
 & & $k=2$ & 0.380 & & 0.190 & \\
 & & $k=3$ & 0.450 &  & 0.225 &  \\
 & & $k=4$ & 0.490 &  & 0.245 & \\
\hline
\end{tabular}
\caption{Comparison of MCM vs ACM thresholds and cost across kernels and profiles. The values of $L^*$ for both method are not displayed for brevity but used to calculate the respective costs.}
\label{tab:table4}
\end{table}
%%%%%%%%%%%%%%%%%%%%%%%%%%%%%%%%%%%%%%%%%%%%%%%%%%%%%%%%%%%%%%
\section{Final remarks and conclusion}\label{sec:conclusion}

This study presents a discrete-space framework for modeling invasion dynamics via  Lattice Difference Equations (LDEs), offering a spatially explicit and computationally tractable alternative to classical reaction-diffusion and integro-difference models. Our results demonstrate that the LDE approach captures essential features of spread—including  threshold-dependent invasion, wave propagation, and infection persistence--while being particularly well-suited to fragmented and patchy landscapes that characterize many real mosquito habitats (\cite{Barton2011BistableWaves, Turelli2017Wolbachia, Riley2003SARS}).

Mathematically, 
the model’s  bi-stability,  driven by the growth function parameters and its induced  Allee effect, yields clear  invasion thresholds: {\it Wolbachia} spreads only when initial infection levels exceed a critical amplitude or spatial extent. These thresholds are consistent with theoretical findings in reaction-diffusion and integro-difference systems (\cite{Li2011Wolbachia, Zhou2008Stability}), but the lattice formulation provides a more ecologically realistic framework in contexts such as urban mosquito habitats or island ecologies (\cite{Yakob2008GMmosquitoes}).
We further show that  traveling wave solutions  exist and propagate at speeds that depend on both the  shape of the dispersal kernel and the  strength of the Allee effect. In particular,  long-tailed dispersal kernels (e.g., Cauchy-type) lead to faster wave speeds  and  lower invasion thresholds compared to exponentially decaying kernels, a pattern supported by prior ecological models (\cite{kot1996dispersal, Clark2001Invasion}). This reinforces the idea that even rare, long-distance dispersal events can play a critical role in accelerating  {\it Wolbachia} spread.

From a public health perspective, our model provides concrete tools for  designing effective release strategies. The  minimum release amplitude and width  required for successful invasion can be computed numerically and used to optimize deployment in heterogeneous environments. Importantly, we extend the model to quantify  local outbreak size  by treating the infection burden at a spatial location as a  sum of Bernoulli-distributed infection events, leading to  Poisson-binomial  distributions. This probabilistic formulation allows us to estimate the likelihood of observing a given infection intensity, which is especially relevant for evaluating field trial outcomes and assessing spatial variability in intervention success (\cite{hancock2011modeling}). This probabilistic  formulation also allows us to propose a novel criterion for release cost minimization
using  modality shift. 

The  discrete nature  of the LDE model enables natural extensions to incorporate  spatial heterogeneity, seasonal forcing, and environmental covariates, such as temperature or rainfall, which are known to modulate  {\it Wolbachia} fitness and transmission (\cite{Murdock2014Temperature}). Moreover, the model may be coupled with  GIS-based data  to simulate invasion in geographically structured regions, enabling  site-specific risk assessments. However, the assumption of uniform dispersal rates (i.e., $\delta_i$ constant)  and symmetric kernels may not reflect true ecological heterogeneity since mosquito movement and release strategies often depend on local environment, barriers, and human behavior. 

Important directions for future research involves considering non constant dispersal rate, non-symmetric kernels, and incorporating  stochasticity  into the LDE framework. Natural mosquito populations are subject to random fluctuations in demography, habitat availability, weather, and human intervention. By extending our model to  stochastic lattice difference equations or  spatially random environments, we can compute  invasion probabilities  rather than deterministic thresholds. This is especially important in ecologically variable regions such as  South Texas, where vector population densities fluctuate seasonally due to rainfall and urban expansion. Incorporating noise also opens the door to modeling stochastic synapse-like rewiring, local extinction-recolonization dynamics, and robustness of wave fronts under uncertainty (\cite{Rohani2010Quality, Barton2011BistableWaves}).

In summary, our LDE-based approach contributes to the growing body of work at the interface of  mathematical biology, spatial ecology, and vector control policy. By providing rigorous analytical results, probabilistic interpretations, and simulation tools tailored to spatially discrete environments, this work helps bridge the gap between theoretical modeling and field-applicable strategy design for {\it Wolbachia}-based interventions against mosquito-borne diseases.

\bibliography{wolbachia_refs}

\section*{Appendix}
\subsection*{Appendix A: Proof of Theorem \ref{thm:persist}}
\begin{proof}
To discuss stability of the LDE, we return to the form

\[
v_i(t+1) = (1 - \delta) f_i(v_i(t)) + \delta \sum_{j \in \mathbb{Z}^d} K_{ij} f_j(v_j(t))\;.
\]
To obtain the Jacobian, we  linearize the system around \( v = v^* \), a fixed point of $f$,  by performing a first-order Taylor expansion:

\[
f(v) \approx f(v^*) + f'(v^*) (v-v^*) + O(v^2).
\]
Since \( v = v^* \) is an equilibrium point (\( f(v^*) = v^* \)), we will have:
\[
f_j(v_j)=f(v_j) \approx f(v_j^*)+f'(v_j^*)( v_j-v^*)=v_j^*+f'(v_j^*) (v_j-v^*)\;.
\]
Substituting this into the LDE equation:
\[
v_i(t+1) \approx (1-\delta)[v_i^*+ f'(v_i^*) (v_i(t)-v^*)] + \delta \sum_{j \in \mathbb{Z}^d} K_{ij}[v_j^*+ f'(v_j^*)( v_j(t)-v^*)].
\]
Put $v=(v_i(t))_{i\in Z^d}$. Rewriting in matrix form, we define the Jacobian matrix \( J \) as:

\[
J_{ij} = \frac{\partial v_i(t+1)}{\partial v_j(t)} \Big|_{v(t)=v^*}.
\]
From the linearized equation, we obtain
\[
J_{ij} = (1-\delta) f'(v_i^*) \delta_{ij} + \delta K_{ij}f'(v_j^*),
\]
where \( \delta_{ij} \) is the Kronecker delta (\( \delta_{ij} = 1 \) if \( i = j \), otherwise 0). We note that 
\begin{itemize}
\item The first term \((1- \delta )f'(v_i^*) \delta_{ij} \) corresponds to  local growth at site \( i \).
\item The second term \( \delta K_{ij} f'(v_j^*) \) accounts for  influence from neighboring sites \( j \) through dispersal.
\end{itemize}
We note that the matrix can be rewritten in operator form as:
\[
J =f'(v^*) \left[(1-\delta)I + \delta K \right] ,
\] 
where:
\begin{itemize}
\item $f'(v^*)=f(v_i^*)_{i\in \Z^d}$.
\item \( I \) is the identity matrix.
\item \( K =(K_{ij})_{i,j\in \Z^d}\) is the dispersal kernel matrix.
\end{itemize}
To find the  eigenvalues $\lambda$ of \( J \),  we write $Jv=\lambda v$, that is $f'(v^*)[(1-\delta)v+\delta Kv]=\lambda v$. Taking the discrete Fourier transform (DFT) of the left and right hand sides and using its linearity, we have that  for some $k\in (\Z^+)^d$,
\[
\widehat{Jv}(k)=f'(v^*)[(1-\delta)\widehat{v}(k)+\delta \widehat{Kv}(k)]=\lambda \widehat{v}(k)\;.
\]
$Kv$ is in fact the convolution between $K$ and $v$. To see this, recall that $K_{ij}=K(x_i-x_j)$ and \[\sum_{j\in \Z^d}K_{ij}v_j=\sum_{j\in \Z^d}K(x_i-x_j)v(x_j)=(K*v)(x_i)\;.\]
It follows that in vector form, we have 
\[
Kv=\left(\sum_{j\in Z^d}K_{ij}v_j\right)_{i\in \Z^d}=((K*v)(x_i))_{i\in \Z^d}=K*v\;.\]
Hence using the property of the Fourier Transform on convolutions, we have:
\[
\widehat{Kv}(k)=\widehat{K}(k)\widehat{v}(k)\;.
\]
In conclusion,  $\widehat{K}(k)$ is an eigenvalue of $K$ and  then \[\lambda_k =  f'(v_k^*) \left[1-\delta+ \delta \widehat{K}(k) \right],\quad k\in (\Z^+)^d\] is an eigenvalue of $J$. For the wave to be stable  stable at $v=v^*$, we require that:
\[
\sup_{k\in (\Z^+)^d} |\lambda_k| <1\;,
\]
which is written as :
\[
\sup_{k\in (\Z^+)^d}  |f'(v_k^*)| |1-\delta + \delta \widehat{K}(k)|  <1\;.
\]

\end{proof}

\subsection*{Appendix B: Proof of Theorem \ref{thm:1}}
\noindent The proof of this Theorem relies on the notion of equicontinuity and  some  classical Theorems, which for sake of completeness will be given below. 
\begin{defn} Let $(X,d_X)$ and $(Y,d_Y)$ be two metric spaces and let $\mathcal{F}$ a collection  of functions from $X$ to $Y$.\
1) $\mathcal{F}$ is point-wise equicontinuous at $x_0\in X$ if 
\[\forall \epsilon>0, \forall f\in \mathcal{F}, \exists \delta=\delta(\epsilon, x_0)>0: \forall x\in X, d_X(x,x_0)\implies d_Y(f(x),f(x_0))<\epsilon\;. \]
2) $\mathcal{F}$ is uniformly equicontinuous if 
\[\forall \epsilon>0, \forall f\in \mathcal{F}, \exists \delta=\delta(\epsilon)>0: \forall x_1,x_2\in X,  d_X(x_1,x_2)\implies d_Y(f(x_1),f(x_2))<\epsilon\;. \]

\end{defn}
\begin{thm}[Arzela-Ascoli Theorem]
Let $X$ be a compact metric space and let $\mathcal{C}(X)$ be the space of continuous functions of $X$. If a sequence $\set{f_n}_{n\in N}$ in $\mathcal{C}(X)$ is bounded and equicontinuous, then it has a uniformly convergence subsequence. 
\end{thm}
\begin{thm}\label{lem1} Let $X$ be a compact metric space. A subset of $\mathcal{C}(X)$ is compact if and only if it closed, bounded and equicontinuous.
\end{thm}
\begin{thm}[Schauder Fixed-Point Theorem] \label{thm:scahuder}Every continuous function from a nonempty convex, compact subset $K$ of a Banach space into itself has a fixed point. 
\end{thm}
\noindent Now were are ready to prove the Theorem.
\begin{proof}
Let $X=[0,1]$, which is a compact metric space of $\R$, as closed bounded subset of $\R$. Let $B=\mathcal{C}(X)$ be the space of continuous functions on $X$. We know that endowed with the norm 
\[\norm{v}_B=\max_{x\in X}|v(x)|\;,\]
$(B, \norm{\cdot}_B)$ is a Banach space.  Now, we consider the subspace $\Gamma=\mathcal{C}_b(X)$ of continuous and uniformly bounded functions on $X$. We define the operator $T: \Gamma\to \Gamma$ as 
\[T(v)=(1-\delta)f(v)+g(v)\;,\] where for any given $v\in \Gamma$, we have 
$g(v): X\to \R$ defined as 
\[g(v)(\xi)=g(v(\xi))=\delta\sum_{j\in \Z^d} K_{\cdot, j}f(v(\xi-j-\cdot))\;.\]
First, we observe that $T$ maps $\Gamma$ into $\Gamma$. Indeed, we recall that $f(0)=0$ and $f(1)=1$. Therefore, 
\[T(0)=(1-\delta_{\cdot})f(0)+g(0)=0+\delta \sum_{j\in \Z^d}K_{\cdot,j}f(0)=0\;.\]
Likewise,
\[T(1)=(1-\delta_{\cdot})f(1)+g(1)=(1-\delta)+\delta\sum_{j\in \Z^d}K_{\cdot,j}f(1)=1-\delta+\delta=1\;.\]
Moreover, for $v\in \Gamma$, we have that for any $\xi\in [0,1]$,  $0\leq v(\xi)\leq 1$ and 
\begin{eqnarray*}
\abs{T(v)(\xi)}&=&\abs{(1-\delta)f(v(\xi))+\delta\sum_{j\in \Z^d} K_{\cdot, j}f(v(\xi-j-\cdot))}\\
&\leq & (1-\delta)+\delta \abs{\sum_{j\in \Z^d} K_{\cdot, j}f(v(\xi-j-\cdot))}, \quad \mbox{since $\abs{f(v(\xi))}<1$}\\
&\leq & (1-\delta)+\delta \sum_{j\in \Z^d} K_{\cdot, j}\abs{f(v(\xi-j-\cdot))}\\
&\leq & (1-\delta)+\delta=1.
\end{eqnarray*}
Therefore, $\norm{T(v)}_B\leq 1$, for all $v\in \Gamma$. Since the functions $v\mapsto f(v)$ and $v\mapsto g(v)$ are continuous, it follows that $T$ is a continuous function from $\Gamma$ to $\Gamma$. Let  $v_1,v_2\in \Gamma$ and let $M$ be uniform bound for elements in $\Gamma$. Consider the convex combination $v(\xi)=(1-\lambda)v_1+\lambda v_2$, for $\lambda\in [0,1]$ and $\xi \in X$. We have that  
\[\abs{v(\xi)}\leq (1-\lambda)\abs{v_1(\xi)}+\lambda \abs{v_2(\xi)}\leq (1-\lambda)M+\lambda M=M\;. \]
This shows that $v\in \Gamma$ and therefore, $\Gamma$ is convex. It remains to show that $\Gamma$ is compact. From Theorem \ref{lem1} above, it is sufficient to show that $\Gamma$ is closed, bounded, and equicontinuous. By definition, $\Gamma$ is closed and bounded. For equicontinuity (which is similar to uniform continuity for functions), we note that given  $\epsilon>0$,  we are looking for $\delta=\delta(\epsilon)$ such $\abs{v(\xi_1)-v(\xi_2)}<\epsilon$, whenever $|\xi_1-\xi_2|<\delta$, and for any $v\in \Gamma$. Since $\Gamma$ is uniformly bounded by a constant $M>0$, and $v\in \Gamma$ is continuous,  therefore, $v$ is uniformly continuous on $X$. Hence, there exists $\delta>0$ such that  $|v(\xi_1)-(\xi_2)|<\epsilon$  whenever $|\xi_1-\xi_2|<\delta$. Theorem \ref{thm:scahuder} allows us to conclude that the function $T$ has a fixed-point, which is a traveling for the {\it Wolbachia}-LDE above.
\end{proof}
 
\subsection*{Appendix C: Proof of Theorem \ref{thm:disOutbreak}}
 
 \begin{proof}
From the definition of $N_i$, we have that
\begin{equation}\label{eqn:minN}
v(t+1,x_i)\neq v(t,x_i)\,\quad \mbox{for $t<N_i$}\;.
\end{equation}
We note that $N_i$ is the time until the system reaches  a fixed point (that is, $v(t,x_i)=0, v(t,x_i)=1$ or $v(t,x_i)=\frac{s_f}{s_h}$) at  location $x_i$, given that is started at $v(0,x_i)$. Since $t\in \Z^+$ takes on discrete values, then $N_i$ is itself distributed according to a geometric distribution. More precisely, calling ``success" the event of reaching a fixed point, given that we started at $v(t,x_i)$, and $q_i=\Pp(\mbox{``success"})=\Pp(\set{v(t+1,x_i)=v(t,x_i)})$, then we have that for $m\in \Z^+, N_i\sim Geom(q_i)$ or that 
\[\Pp(N_i=m)=(1-q_i)^mq_i\;.\]
 Define 
\[p_{it}=v(t,x_i)\;.\]
 We know that if $p_{i0}=p_{i1}=\cdots=p_{i,N_i-1}=p$, then $Y_i|N_i\sim Binom(N_i,p)$, but from \eqref{eqn:minN}, this case is impossible. That leaves the general case where $0\neq p_{i0}\neq p_{i1}\neq \cdots\neq p_{i,N_i-1}\neq 1$ and the particular case where the $p_{it}$'s are  small and $N_i$ is large. In the former case, as the sum of independent Bernoulli distributions with non  constant probability of success $p_{it}$, we have that  
 \[
Y_i|N_i=\sum_{t=0}^{N_i-1}Y_{it}\sim PB(p_{i0}, p_{i1},\cdots, p_{i,N_i-1})\;.
\]
 In the latter case, $Y_i|N_i\sim Poisson (\lambda_i)$, where $\lambda_i=\sum_{t=0}^{N_i-1}p_{it}$ is  referred to as the infection burden at location $x_i$. More precisely, using the law of total probability, we have that 
\begin{equation}\label{eqn:distOfOutbreal}
\Pp(Y_i=k)=\sum_{m=0}^{\infty}\Pp(Y_i=k|N_i=m)\Pp(N_i=m)
\end{equation}
And the first term in the above equality is obtained using the classic definition of the Poisson-Binomial distribution as: 
 \begin{equation*}
 \Pp(Y_i=k|N_i=m)=\sum_{U\in \Lambda_k}\prod_{j\in U}p_{ij}\prod_{i\in U^c}(1-p_{ij}) \;,
 \end{equation*}
 where 
 \[\Lambda_k=\set{S_k:  S_k\subseteq \set{1,2,\cdots, m}, ~|S_k|=k}\;,\] where  $|A|$ is the cardinality of the set $A$.
 \end{proof}

\end{document}